\documentclass[11pt,reqno]{amsart}
\usepackage{amsfonts}
\usepackage{mathrsfs}
\usepackage{amsmath}
\usepackage{amssymb,amsfonts}
\usepackage{amsthm}
\usepackage{graphicx}
\usepackage{extarrows}
\usepackage{mathtools}
\usepackage{enumerate}
\usepackage{appendix}

\newtheorem{thm}{Theorem}[section]
\newtheorem{lem}{Lemma}[section]
\newtheorem{defn}{Definition}[section]
\newtheorem{prop}{Proposition}[section]
\newtheorem{coro}{Corollary}[section]\numberwithin{equation}{section}
\newtheorem{rmk}{Remark}[section]

\def\pf{{\textit {Proof:} }}

\newcommand{\mysection}[1]{\section{#1}\setcounter{equation}{0}}

\newfont{\bb}{msbm10 at 11pt}

\newcommand{\bal}{\begin{aligned}}      \newcommand{\eal}{\end{aligned}}
\newcommand{\ba}{\begin{array}}      \newcommand{\ea}{\end{array}}
\newcommand{\bc}{\begin{center}}     \newcommand{\ec}{\end{center}}
\newcommand{\be}{\begin{enumerate}}  \newcommand{\ee}{\end{enumerate}}
\newcommand{\beq}{\begin{eqnarray}}  \newcommand{\eeq}{\end{eqnarray}}
\newcommand{\beQ}{\begin{eqnarray*}} \newcommand{\eeQ}{\end{eqnarray*}}
\newcommand{\bi}{\begin{itemize}}    \newcommand{\ei}{\end{itemize}}
\newcommand{\bt}{\begin{tabular}}    \newcommand{\et}{\end{tabular}}
\newcommand{\bdm}{\begin{displaymath}} \newcommand{\edm}{\end{displaymath}}

\newcommand{\lrw}{\longrightarrow}

\newcommand{\Llrw}{\Longleftrightarrow}

\def\qed{\hfill{Q.E.D.}\smallskip}

\newcommand{\ls}{\setlength{\baselineskip}{12pt}
                 \setlength{\parskip}{3mm}}

\begin{document}

\allowdisplaybreaks

\title[Intrinsic NC geometry]{Deformation quantization and intrinsic noncommutative differential geometry}


\author[H Gao]{Haoyuan Gao$^{3,1,2}$}
\author[X Zhang]{Xiao Zhang$^{1,2,4}$}

\address[]{$^{1}$Academy of Mathematics and Systems Science, Chinese Academy of Sciences, Beijing 100190, PR China}
\address[]{$^{2}$School of Mathematical Sciences, University of Chinese Academy of Sciences, Beijing 100049, PR China}
\address[]{$^{3}$Shanghai Center for Mathematical Sciences, Fudan University, Jiangwan Campus, Shanghai 200438, PR China}
\address[]{$^{4}$Guangxi Center for Mathematical Research, Guangxi University, Nanning, Guangxi 530004, PR China}

\email{hy\_gao@fudan.edu.cn$^{3}$, gaohaoyuan@amss.ac.cn$^{1,2}$}
\email{xzhang@amss.ac.cn$^{1,2}$, xzhang@gxu.edu.cn$^{4}$}

\date{}

\begin{abstract}

We provide an intrinsic formulation of the noncommutative differential geometry developed earlier by Chaichian, Tureanu, R. B. Zhang and the second author. This yields geometric definitions of covariant derivatives of noncommutative metrics and curvatures, as well as the noncommutative version of the first and the second Bianchi identities. Moreover, if a noncommutative metric and chiral coefficients satisfy certain conditions which hold automatically for quantum fluctuations given by isometric embedding, we prove that the two noncommutative Ricci curvatures are essentially equivalent. For (pseudo-) Riemannian metrics given by certain type of spherically symmetric isometric embedding, we compute their quantum fluctuations and curvatures. We find that they have closed forms, which indicates that the quantization of gravity is renormalizable in this case. Finally, we define quasi-connections and their curvatures with respect to general associative star products constructed by Kontsevich on Poisson manifolds. As these star products are not compatible with the Leibniz rule, we can only prove the first Bianchi identity.

\end{abstract}

\maketitle \pagenumbering{arabic}

\mysection{Introduction}
\ls

Gravity is essentially a theory of spacetime geometry. In the concept of quantum effects of gravity, the Heisenberg uncertainty relations would result in noncommutativity of spacetime variables for sufficiently small distances. In 1947, Snyder, C.N. Yang made the first attempts to quantize spacetimes \cite{S, Y}, which are referred as Snyder's quantum space-times and Yang's quantum phase spaces \cite{LW1, LW2}. In their approach, spacetime variables were represented by Hermitian operators with discrete eigenvalues. This idea to encode geometry of a space by its algebras of functions was realized prominently by Connes to establish noncommutative geometry using spectral triples \cite{C}, where the main ingredients are the noncommutative analog of the Dirac operator acting on a representation space of the algebra, the spectrum of this generalized Dirac operator and the cyclic (co)homology. They are used to encode the information of noncommutative manifold structure, noncommutative metric and noncommutative curvature respectively. The overview of its applications to physics can be found in \cite{CCS}.

However, the metric and curvature information in an infinitesimal neighborhood of manifold is still lack as it is not known what means to take derivatives when coordinate variables are operators. Alternatively, deformation quantization deforms the commutative algebras of functions based on pointwise commutative multiplication to noncommutative algebras of functions based on certain noncommutative products such as the Moyal product, but still keeps spacetime variables usual functions, c.f. \cite{BFFLS1, BFFLS2}. In recent years, there have been intensive research activities on noncommutative gravity in frame of deformation quantization, c.f. \cite{SW, ABDMSW, ADMW} and references therein, where general relativity is adopted to the noncommutative setting in an intuitive way, as pointed out in \cite{TZ}.

In \cite{CTZZ, WZZ, WZZ1, ZZ}, a mathematically rigorous and complete theory of noncommutative differential geometry was developed on a coordinate chart $U$ of a (pseudo-) Riemannian manifold. The idea is to embed $U$ isometrically into a flat (pseudo-) Euclidean space and use the isometric embedding to construct the noncommutative analogues of metric, connection and curvature. They yield the noncommutative Einstein field equations. It was found that the deformation quantization of the Schwarzschild metric does not depend on time and yields an unevaporated quantum black hole \cite{WZZ}, and the quantum fluctuation of the plane-fronted gravitational wave is the exact solution of the noncommutative vacuum Einstein field equations \cite{WZZ1}. We refer to \cite{PT, R} for the review on general existence of isometric embedding and applications in physics.

In this paper, we provide the intrinsic theory of noncommutative differential geometry via deformation quantization, without using the isometric embedding. In Section 2, we review basic knowledge on noncommutative metrics, connections, curvatures and state the main theorem. In Section 3, we study the intrinsic formulation of covariant derivatives of noncommutative metrics and curvatures from the geometric point of view. This yields noncommutative version of the first and the second Bianchi identities. In Section 4, we show the two noncommutative Ricci curvatures are essentially equivalent if noncommutative metrics and chiral coefficients satisfy certain conditions. These conditions hold automatically for quantum fluctuations given by isometric embedding. In Section 5, we show that the quantum fluctuations and their curvatures have closed forms coming from Moyal products of trigonometric functions if (pseudo-) Riemannian metrics are given by certain type of spherically symmetric isometric embedding. In Section 6, we define quasi-connections and their curvatures with respect to general associative star products and prove the first Bianchi identity.

\mysection{Basic knowledge and main theorem}

In this section, we review basic knowledge on the noncommutative differential geometry and, in particular, recall the intrinsic setting of noncommutative differential geometry proposed by the second author \cite{Z}, without using the isometric embedding. The theory is stated in the spirit of those in \cite{CTZZ, WZZ, WZZ1, ZZ, Z} with necessary modification and improvement.

Let $M$ be an $n$-dimensional differentiable manifold and $U \subset M$ be a coordinate chart equipped with natural coordinates $(x^1, \cdots, x^n)$. Let $\hbar$ be the Planck constant viewed as an indeterminate. Denote $\mathbb{R}[[\hbar]]$ the ring of formal power series in $\hbar$ with real coefficients, and $\mathcal{A}_U$ the set of formal power series in $\hbar$ with coefficients being real smooth functions on $U$
$$
\mathcal{A}_U = C^{\infty}(U)[[\hbar]]=\Big\{ \sum\limits_{q = 0}^{\infty}f_q\hbar^q \Big| f_q \in C^{\infty}(U)   \Big\}.
$$
$\mathcal{A}_U$ is an $\mathbb{R}[[\hbar]]$-module.

Throughout the paper, all the indices $i$, $j$, $k$, $l$, $\cdots$, range from $1$ to $n$, $q \in \mathbb{N}_0$. We also use the Einstein summation convention. Given two smooth functions $u$, $v$ on $U$, we denote $uv$ their usual pointwise product. For any skew-symmetric $n \times n$ real constant matrix $(\theta^{ij})$ on $U$, the Moyal product of $u$ and $v$ with respect to $(\theta^{ij})$ is defined as
\begin{align}
    (u * v)(x) = \Big[\exp(\hbar\theta^{ij}\partial_i\partial'_j)u(x)v(x')\Big]_{x = x'},
\end{align}
where $x$ and $x'$ denote the same coordinate system and
\begin{align*}
\partial_i = \frac{\partial}{\partial{x^i}}, \quad \partial'_i = \frac{\partial}{\partial{(x')^i}}.
\end{align*}
It is clearly that
\begin{align*}
u*v \in \mathcal{A}_U.
\end{align*}

Extending by $\mathbb{R}[[\hbar]]$-bilinearity, the Moyal product provides an associative $\mathbb{R}[[\hbar]]$-bilinear product on $\mathcal{A}_U$, c.f. \cite{LPV}. The Moyal algebra is $\mathcal{A}_U$ equipped with the Moyal product, which is a formal deformation of the algebra of real smooth functions on $U$.

Extend $\partial_i$ to $\mathcal{A}_U$ by $\mathbb{R}[[\hbar]]$-linearity, the Moyal product satisfies
\bi
    \item[(i)] Noncommutativity: $[x^i, x^j]=x^i * x^j - x^j * x^i = 2\hbar\theta^{ij}$;\\
    \item[(ii)] Leibniz rule: $\partial_i(u * v) = (\partial_i{u}) * v + u * (\partial_i{v})$, for $u, v \in \mathcal{A}_U$.
\ei

The following definitions of noncommutative tangent bundles and metrics can be found, e.g., in \cite{CTZZ, Z}. Denote
\begin{align*}
E_i = \tilde{E}_i =\partial_i, \quad 1 \leq i \leq n.
\end{align*}

\begin{defn}
The noncommutative left ({\rm{resp.}} right) tangent bundle $\mathcal{T}_U$ ({\rm{resp.}} $\tilde{\mathcal{T}}_U$) on $U$ is the free left ({\rm{resp.}} right) $\mathcal{A}_U$-module with basis $\{E_1, \cdots, E_n\}$ ({\rm{resp.}} $\{\tilde{E}_1, \cdots, \tilde{E}_n\}$), i.e.,
\begin{align*}
\mathcal{T}_U & = \Big\{a^i * E_i \,\Big | \,a^i \in \mathcal{A}_U, \, a^i * E_i = 0 \Longleftrightarrow a^i = 0. \Big\},\\
\tilde{\mathcal{T}}_U & = \Big\{\tilde{E}_i * \tilde{a}^i \,\Big | \, \tilde{a}^i \in \mathcal{A}_U, \, \tilde{E}_i * \tilde{a}^i = 0 \Longleftrightarrow \tilde{a}^i = 0. \Big\}.
\end{align*}
An element of $\mathcal{T}_U$ ({\rm{resp.}} $\tilde{\mathcal{T}}_U$) is called a left ({\rm{resp.}} right) vector field.
\end{defn}

\begin{defn}
A noncommutative metric $g$ on $U$ is a homomorphism of two-sided $\mathcal{A}_U$-modules
\begin{align*}
g : \mathcal{T}_U \otimes_{\mathbb{R}[[\hbar]]} \tilde{\mathcal{T}}_U \longrightarrow \mathcal{A}_U,
\end{align*}
such that the matrix
\begin{align*}
(g_{ij}) \in \mathcal{A}_U ^{n \times n}, \quad g_{ij} = g(E_i, \tilde{E}_j)
\end{align*}
is invertible, i.e., there exists a unique matrix $(g^{ij}) \in \mathcal{A}_U ^{n \times n}$ such that
\begin{align*}
g_{ik} * g^{kj} = g^{jk} * g_{ki} = \delta_i^j.
\end{align*}
\end{defn}

Let $(g_l^{ij})$ be the left inverse of $(g_{ij})$ and $(g_r^{ij})$ be the right inverse of $(g_{ij})$. Since the Moyal product is associative,
\begin{align*}
    g_l^{ij} = g_l^{ip} * \delta_p^j = g_l^{ip} * g_{pk} * g_r^{kj} = \delta_k^i * g_r^{kj} = g_r^{ij}.
\end{align*}
Therefore the left inverse and the right inverse coincide.

In classical differential geometry, the cotangent bundle is the dual of the tangent bundle. Inspired by this, we can define the noncommutative cotangent bundles as the dual modules of the noncommutative tangent bundles. As the dual of a left (\rm{resp.} right) $\mathcal{A}_U$-module is a right (\rm{resp.} left) $\mathcal{A}_U$-module and the dual of a free module is also free, we may use the noncommutative metric $g$ to induce bases of the cotangent bundles dual to $E_i$ and $\tilde{E}_j$ respectively, i.e., let $E^i$, $\tilde{E}^j$ be dual bases of $\tilde{E}_j$, $E_i$ respectively, we have
\begin{align*}
    g(E^i, \tilde{E}_j) = g(E_j, \tilde{E}^i)= \delta_j^i.
\end{align*}

\begin{defn}
The noncommutative left ({\rm{resp.}} right) cotangent bundle $\mathcal{T}_U^*$ ({\rm{resp.}} $\tilde{\mathcal{T}}_U^*$) on $U$ with respect to the noncommutative metric $g$ is the free left ({\rm{resp.}} right) $\mathcal{A}_U$-module with basis $\{E^1, \cdots, E^n\}$ ({\rm{resp.}} $\{\tilde{E}^1, \cdots, \tilde{E}^n\}$)
\begin{align*}
\mathcal{T}_U^* & = \Big\{a_i * E^i \,\Big | \,a_i \in \mathcal{A}_U, \, a_i * E^i = 0 \Longleftrightarrow a_i = 0 \Big\},\\
\tilde{\mathcal{T}}_U^* & = \Big\{\tilde{E}^i * \tilde{a}_i \,\Big | \, \tilde{a}_i \in \mathcal{A}_U, \, \tilde{E}^i * \tilde{a}_i = 0 \Longleftrightarrow \tilde{a}_i = 0 \Big\}.
\end{align*}
\end{defn}

The left ({\rm{resp.}} right) cotangent bundle is the dual of the right ({\rm{resp.}} left) tangent bundle. Analogous to the classical situation, the noncommutative metric $g$ acts as an element of $\tilde{\mathcal{T}}_U^* \otimes_{\mathcal{A}_U} \mathcal{T}_U^*$,
\begin{align}
    \tilde{E}^i \otimes g_{ij} * E^j = \tilde{E}^i * g_{ij} \otimes E^j.
\end{align}
The inverse matrix $(g^{ij})$ can be viewed as a homomorphism of two-sided modules
\begin{align*}
    g^{-1} : \mathcal{T}_U^* \otimes_{\mathbb{R}[[\hbar]]} \tilde{\mathcal{T}}_U^* \longrightarrow \mathcal{A}_U
\end{align*}
such that
\begin{align*}
g^{-1}(E^i, \tilde{E}^j) = g^{ij}.
\end{align*}
Similarly, $g^{-1}$ acts as an element of $\tilde{\mathcal{T}}_U \otimes_{\mathcal{A}_U} \mathcal{T}_U$,
\begin{align*}
    \tilde{E}_i \otimes g^{ij} * E_j = \tilde{E}_i * g^{ij} \otimes E_j.
\end{align*}

\begin{defn}
A noncommutative left ({\rm{resp.}} right) connection $\nabla $ is a map
\begin{align*}
\nabla: \mathcal{T}_U \longrightarrow \tilde{\mathcal{T}} ^* _U \otimes _{\mathcal{A}_U} \mathcal{T}_U \quad
({\rm{resp.}} \,\,\tilde{\nabla}: \tilde{\mathcal{T}} _U \longrightarrow \tilde{\mathcal{T}} _U \otimes _{\mathcal{A}_U} \mathcal{T} ^* _U)
\end{align*}
such that noncommutative left ({\rm resp.} right) covariant derivatives
\begin{align*}
\nabla_i : \mathcal{T}_U \longrightarrow \mathcal{T}_U, \quad
({\rm{resp.}}\,\,\tilde{\nabla}_i: \tilde{\mathcal{T}} _U \longrightarrow \tilde{\mathcal{T}} _U);
\end{align*}
defined by
\begin{align*}
\nabla_i V =g(E_i, \tilde{E} ^k) * W_k \quad ({\rm{resp.}} \,\,\tilde{\nabla}_i \tilde{V} =\tilde{W} _k * g(E^k, \tilde{E}_i))
\end{align*}
for any
\begin{align*}
\nabla V =\tilde{E} ^k \otimes W_k,\,\, W_k \in \mathcal{T}_U \quad  ({\rm{resp.}} \,\,\tilde{\nabla} \tilde{V} =\tilde{W} _k \otimes E^k,\,\, \tilde{W}_k \in \tilde{\mathcal{T}}_U)
\end{align*}
satisfy
\bi
\item[(i)] $\mathbb{R}[[\hbar]]$-linearity: For $a, \,b \in \mathbb{R}[[\hbar]]$, $V, \,W \in \mathcal{T}_U$ ({\rm{resp.}} $\tilde{V}, \,\tilde{W} \in \tilde{\mathcal{T}}_U$),
\begin{align*}
        \nabla_i(aV + bW)  & = a\nabla_iV + b\nabla_i W,\\
      \tilde{\nabla} _i(a \tilde{V} + b \tilde{W}) & = a \tilde{\nabla}_i \tilde{V} + b \tilde{\nabla}_i \tilde{W};
\end{align*}
\item[(ii)] Leibniz rule: For $f \in \mathcal{A}_U$, $V \in \mathcal{T}_U$ ({\rm{resp.}} $\tilde{V} \in \tilde{\mathcal{T}}_U$),
\begin{align*}
        \nabla_i(f * V) & = (\partial_i f) * V + f * \nabla_iV,\\
       \tilde{\nabla} _i(\tilde{V} * f) & = \tilde{V} *(\partial_i f) + \tilde{\nabla} _i \tilde{V} *f.
\end{align*}
\ei
\end{defn}

\begin{rmk}
The noncommutative left ({\rm{resp.}} right) covariant derivative along a left ({\rm{resp.}} right) vector field $V = a^i * E_i$ ({\rm{resp.}} $\tilde{V} = \tilde{E}_i * \tilde{a}^i$) with $a^i \in \mathcal{A}_U$ (\rm{resp.}\,\,$\tilde{a}^i \in \mathcal{A}_U$) is defined as the $\mathbb{R}[[\hbar]]$-linear map
\begin{align*}
    \nabla_V : \mathcal{T}_U \to \mathcal{T}_U \quad ({\rm{resp.}} \quad \it{\tilde{\nabla}_{\tilde{V}}: \tilde{\mathcal{T}}_U \to \tilde{\mathcal{T}}_U})
\end{align*}
given by
\begin{align*}
    \nabla_V X = a^i * (\nabla_iX) \quad ({\rm{resp.}} \quad \it{\tilde{\nabla}_{\tilde{V}} \tilde{X} = (\tilde{\nabla}_i \tilde{X}) * \tilde{a} ^i})
\end{align*}
for $X \in \mathcal{T}_U$ ({\rm resp.} $\tilde{X} \in \tilde{\mathcal{T}}_U$). But they are not compatible with the Leibniz rule. Indeed, let $X=b^j * E_j$,  $b^j \in \mathcal{A}_U$. By the Leibniz rule,
\begin{align*}
  \nabla_V(b^j * E_j)=&(V b^j) * E_j + b^j * \nabla_V E_j \\
                           =& a^i * \partial_i b^j * E_j + b^j * a^i * \nabla_{i} E_j.
\end{align*}
On the other hand, by the definition,
\begin{align*}
    \nabla_V(b^j * E_j) = & a^i * \nabla_i (b^j * E_j)\\
                        = & a^i * \partial_i b^j * E_j + a^i * b^j * \nabla_{i} E_j.
\end{align*}
They are not equal to each other unless
\begin{align*}
a^i * b^j = b^j * a^i
\end{align*}
for all $i$, $j$. This is generally impossible. As a consequence, it indicates that the noncommutative covariant derivatives are not well-defined with respect to orthonormal basis.
\end{rmk}

The left and right connections are uniquely determined by connection coefficients $\Gamma_{ij}^k$ and $\tilde{\Gamma}_{ij}^k$, which are elements of $\mathcal{A}_U$
\begin{align*}
\nabla_i E_j = \Gamma_{ij}^k * E_k, \quad \tilde{\nabla}_i\tilde{E}_j = \tilde{E}_k * \tilde{\Gamma}_{ij}^k.
\end{align*}

Similar to the classical differential geometry, a noncommutative left connection $\nabla $ on the left tangent bundle induces a unique noncommutative right connection $\tilde{\nabla}$ on the right cotangent bundle $\tilde{\mathcal{T}}_U^*$ in terms of noncommutative metric $g$
\begin{align*}
    \partial_i g(E_k, \tilde{E}^j) = g(\nabla_i E_k, \tilde{E}^j) + g(E_k, \tilde{\nabla}_i \tilde{E}^j).
\end{align*}
It yields
\begin{align*}
    \tilde{\nabla} _i\tilde{E}^j = -\tilde{E}^k * \Gamma_{ik}^j.
\end{align*}
Moreover, a noncommutative right connection $\tilde{\nabla}$ on the right tangent bundle also induces a noncommutative left connection $\nabla$ on the left cotangent bundle which yields
\begin{align*}
    \nabla_iE^j = -\tilde{\Gamma}_{ik}^j * E^k.
\end{align*}

Inspired by the Levi-Civita connection of a (pseudo-) Riemannian metric, the second author introduced the canonical connection \cite{Z}.

\begin{defn}
Given a noncommutative metric $g$ and a set of elements $\Upsilon_{ijk}$ of $\mathcal{A}_U$ with $$\Upsilon_{ijk} = \Upsilon_{jik},$$which are referred as {\em chiral coefficients}. Denote
\begin{align*}
    \Gamma_{ijk} = \Gamma_{ij}^l * g_{lk}, \quad \tilde{\Gamma}_{ijk} = g_{kl} * \tilde{\Gamma}_{ij}^l.
\end{align*}
A noncommutative connection, which consists of a noncommutative left connection $\nabla $ and a noncommutative right connection $\tilde{\nabla}$, is canonical with respect to $g$ and $\Upsilon_{ijk}$ if it satisfies
\bi
\item[(i)] Compatibility: $ \partial_k g_{ij} = g(\nabla_k E_i, \tilde{E}_j) + g(E_i, \tilde{\nabla} _k \tilde{E} _j) = \Gamma_{kij} + \tilde{\Gamma}_{kji};$
\item[(ii)] Torsion free: $ \nabla_i E_j = \nabla_j E_i $,  $\tilde{\nabla}_i\tilde{E}_j = \tilde{\nabla}_j\tilde{E}_i;$
\item[(iii)] Chirality: $ \Gamma_{ijk} - \tilde{\Gamma}_{ijk} = \Upsilon_{ijk}.$
\ei
\end{defn}

The torsion free condition implies
\begin{align*}
    \Gamma_{ij}^k = \Gamma_{ji}^k, \quad \tilde{\Gamma}_{ij}^k = \tilde{\Gamma}_{ji}^k.
\end{align*}
It is straightforward that
\beq\label{Gamma-1}
\begin{aligned}
   2 \Gamma_{ijk} & = \partial_i g_{jk} + \partial_j g_{ki} - \partial_kg_{ij} + \Upsilon_{ikj} + \Upsilon_{jik} - \Upsilon_{kji}\\
                  & = \partial_i g_{jk} + \partial_j g_{ki} - \partial_kg_{ji} + \Upsilon_{ijk}\\
                  & = \partial_i\Big(\frac{g_{jk} + g_{kj}}{2}\Big) + \partial_j\Big(\frac{g_{ki} + g_{ik}}{2}\Big)
                      - \partial_k\Big(\frac{g_{ij} + g_{ji}}{2}\Big) +\Upsilon_{ijk},
\end{aligned}
\eeq
and
\beq\label{Gamma-2}
\begin{aligned}
2 \tilde{\Gamma}_{ijk} & = \partial_i g_{jk} + \partial_jg_{ki} - \partial_kg_{ij} + \Upsilon_{ikj} - \Upsilon_{jik} - \Upsilon_{kji}\\
                       & = \partial_i g_{jk} + \partial_j g_{ki} - \partial_kg_{ji} -\Upsilon_{ijk}\\
                       & = \partial_i\Big(\frac{g_{jk} + g_{kj}}{2}\Big) + \partial_j\Big(\frac{g_{ki} + g_{ik}}{2}\Big)
                           - \partial_k\Big(\frac{g_{ij} + g_{ji}}{2}\Big) - \Upsilon_{ijk}.
\end{aligned}
\eeq
In classical Riemannian geometry, the chiral coefficients vanish and $\Gamma_{ijk}$ reduce to the Christoffel symbols.

For any $f \in \mathcal{A}_U$, it is easy to verify
\begin{align*}
[E_i, E_j] f= [\tilde{E}_i, \tilde{E} _j]f =\partial_i \partial _j f -\partial _j \partial _i f =0.
\end{align*}
Thus the left curvature operators $\mathcal{R}_{E_iE_j}$ and the right curvature operators $\tilde{\mathcal{R}}_{\tilde{E}_i\tilde{E}_j}$ can be defined as the following $\mathcal{A}_U$-linear operators
\begin{align*}
    \mathcal{R}_{E_iE_j} &= [\nabla_i, \nabla_j]:\,\, \mathcal{T}_U \lrw \mathcal{T}_U,\\
    \tilde{\mathcal{R}}_{\tilde{E}_i\tilde{E}_j} &= [\tilde{\nabla}_i, \tilde{\nabla}_j] : \,\,\tilde{\mathcal{T}}_U \lrw \tilde{\mathcal{T}}_U.
\end{align*}
For the canonical connection, the left Riemannian curvatures $R_{lkij}$ and right Riemannian curvatures $\tilde{R}_{lkij}$ are defined as
\begin{align*}
    R_{lkij} = g(\mathcal{R}_{E_iE_j}E_k, \tilde{E}_l), \quad
    \tilde{R}_{lkij} = -g(E_k, \tilde{R}_{\tilde{E}_i\tilde{E}_j}\tilde{E}_l).
\end{align*}
They satisfy
\begin{align*}
    R_{lkij} = -R_{lkji}=\tilde{R}_{lkij}, \quad   R_{lkij} \not \equiv - R_{klij}.
\end{align*}
Therefore the left curvatures are sufficient for the purpose. There are two Ricci curvatures $R_{kj}$ and $\Theta_{il}$ obtained by contracting $l$, $i$ and $k$, $j$ in $R_{lkij}$ respectively
\begin{align*}
    R_{kj} = &g(\mathcal{R}_{E_iE_j}E_k, \tilde{E}_l) * g^{li} = R_{lkij} * g^{li},\\
    \Theta_{il} = &g^{jk} * g(\mathcal{R}_{E_iE_j}E_k, \tilde{E}_l) = g^{jk} * R_{lkij}.
\end{align*}
Raising the index at $k$ and $l$ respectively, we have Ricci curvatures
\begin{align*}
    R_j^p = &g^{pk} * g(\mathcal{R}_{E_iE_j}E_k, \tilde{E}_l) * g^{li}= g^{pk} * R_{lkij} * g^{li},\\
    \Theta_i^p = & g^{jk} * g(\mathcal{R}_{E_iE_j}E_k, \tilde{E}_l) * g^{lp} = g^{jk} * R_{lkij} * g^{lp}.
\end{align*}
The two Ricci curvatures $R_i^p$ and $\Theta_i^p$ are not equal to each other in the noncommutative case. But their traces coincide and yield the same scalar curvature
\begin{align*}
    R = R_j ^j = \Theta_i ^i.
\end{align*}

As elements of $\mathcal{A}_U$, there are the following power series expansions
\begin{align}
g_{ij} &= \sum\limits_{q = 0}^{\infty}g_{ij}[q]\hbar^{q}, \quad \quad \,\,\,\,g_{ij}[q] \in C^{\infty}(U), \label{gij-q1} \\
\Upsilon _{ijk}&=\sum\limits_{q = 0}^{\infty}\Upsilon _{ijk}[q]\hbar^{q}, \quad \, \,\Upsilon _{ijk}[q] \in C^{\infty}(U), \label{Uijk}\\
R_{lkij} &= \sum\limits_{q = 0}^{\infty}R_{lkij}[q]\hbar^q, \quad R_{lkij}[q] \in C^{\infty}(U), \label{riem}\\
R_{ij} &= \sum\limits_{q = 0}^{\infty}R_{ij}[q]\hbar^q, \quad \quad \,\,\,R_{ij}[q] \in C^{\infty}(U),   \label{lowricci-q1}\\
\Theta_{ij} &= \sum\limits_{q = 0}^{\infty}\Theta_{ij}[q]\hbar^q, \quad \quad \,\,\,\Theta_{ij}[q] \in C^{\infty}(U), \label{lowricci-q2}\\
R^i_j &= \sum\limits_{q = 0}^{\infty}R^i_j[q]\hbar^q, \quad \qquad  R^i_j[q] \in C^{\infty}(U),   \label{ricci-q1}\\
\Theta^i_j &= \sum\limits_{q = 0}^{\infty}\Theta^i_j[q]\hbar^q, \quad \qquad \Theta^i_j[q] \in C^{\infty}(U). \label{ricci-q2}
\end{align}

In this paper, we prove the following theorem.

\begin{thm}
Let $M$ be an $n$-dimensional smooth manifold and $U \subset M$ a coordinate chart. Let $\nabla $, $\tilde{\nabla}$ be the canonical connections with respect to noncommutative metric $g$ and chiral coefficients $\Upsilon_{ijk}$ on $U$. If $g_{ij}$ satisfy
\begin{align}
    g_{ij}[2q] = g_{ji}[2q], \quad g_{ij}[2q + 1] =  -g_{ji}[2q + 1], \label{ricci-g}
\end{align}
and $\Upsilon_{ijk}$ satisfy
\begin{align}
 \Upsilon_{ijk}[2q] = 0,   \label{ricci-U}
\end{align}
then two Ricci curvatures are equivalent in the sense that
\begin{align}
    R_{ij}[2q] = \Theta_{ji}[2q], \quad R_{ij}[2q + 1] = - \Theta_{ji}[2q + 1]
\end{align}
and
\begin{align}
    R^i_j[2q] = \Theta^i_j[2q], \quad \,\,\, R^i_j[2q + 1] = -\Theta^i_j[2q + 1].
\end{align}
\end{thm}

In particular, if noncommutative metric and chiral coefficients are given by an isometric embedding, then (\ref{ricci-g}), (\ref{ricci-U}) hold and the theorem follows.

Finally, we would like to remark that, in Poisson geometry, the Moyal product is a deformation quantization of the constant Poisson structure $$\pi =\frac{1}{2}\theta^{ij}\partial_i \wedge \partial_j$$ for constant skew-symmetric matrix $(\theta^{ij})$. If $\theta^{ij}$ are smooth functions, $\pi$ still gives a Poisson structure if its Schouten-Nijenhuis bracket vanishes, $$[\pi, \pi]_S = 0.$$
However, the corresponding Moyal product is not associative. In the pioneer work, Kontsevich proved that there always exists an associative noncommutative star product which provides the deformation quantization for any Poisson structure \cite{Ko}. A concrete algorithm for calculating those integrals appeared in Kontsevich's formula was given by Banks, Panzer and Pym using integer-linear combinations of multiple zeta values \cite{BPP}, which produce the first software package for the symbolic calculation of Kontsevich's formula. It should be pointed out that Kontsevich's star products are not compatible with the Leibniz rule. It indicates the theory of noncommutative differential geometry depends on the choice of coordinate systems in $U$. As coordinate systems correspond to observers, this fits Bohr's opinion that evidence obtained under different experimental conditions cannot be comprehended within a single picture, but must be regarded as complementary in the sense that only the totality of the phenomena exhausts
the possible information about the objects.

\mysection{Curvature operators and Bianchi identities}
\ls

In this section, we study the covariant derivatives of noncommutative metrics and curvatures from the geometric point of view. This yields noncommutative version of the first and the second Bianchi identities.

\begin{prop}\label{gij0}
Let $M$ be an $n$-dimensional differentiable manifold and $U \subset M$ be a coordinate chart equipped with natural coordinates $(x^1, \cdots, x^n)$. Let $g$ be a homomorphism of two-sided $\mathcal{A}_U$-modules given by (\ref{gij-q1}) where $(g_{ij}[0])$ is not necessarily symmetric. If $(g_{ij}[0])$ is invertible on $U$ with the inverse matrix $(g^{ij}[0])$, then (\ref{gij-q1}) gives a noncommutative metric $g$ on $U$.
\end{prop}
\pf For any smooth functions $u(x)$, $v(x)$ over $U$, denote
\begin{align}
    \mu_q (u, v)(x) = \frac{1}{q!}\Big[(\theta^{ij}\partial_i\partial'_j)^q u(x)v(x')\Big]_{x = x'}.    \label{mun}
\end{align}
Let $g^{ij}$ have the power series expansions
\begin{align}
g^{ij} &= \sum\limits_{q = 0}^{\infty}g^{ij}[q]\hbar^{q} \in \mathcal{A}_U, \quad g^{ij}[q] \in C^{\infty}(U).\label{gij-q2}
\end{align}
By viewing $g^{ij}$ as the right inverse, we obtain the recursive formula for $q \in \mathbb{N}$
\beq
\begin{aligned}
g^{ij}[q] = &-\sum\limits_{r = 1}^qg^{ik}[0]g_{kl}[r]g^{lj}[q - r]\\
                     &- \sum\limits_{r = 1}^q\sum\limits_{s = 0}^{q - r}g^{ik}[0]\mu_r \Big(g_{kl}[s], g^{lj}[q - r - s] \Big).  \label{r-g}
\end{aligned}
\eeq
On the other hand, by viewing $g^{ij}$ as the left inverse, we obtain
\beq
\begin{aligned}
g^{ij}[q] = &- \sum\limits_{r = 1}^qg^{ik}[q - r]g_{kl}[r]g^{lj}[0] \\
                &- \sum\limits_{r = 1}^q\sum\limits_{s = 0}^{q - r}\mu_r \Big(g^{ik}[q - r - s], g_{kl}[s] \Big)g^{lj}[0]. \label{l-g}
\end{aligned}
\eeq
Thus the matrix $(g_{ij})$ is invertible in $\mathcal{A}_U ^{n \times n}$ if and only if the matrix $(g_{ij}[0](x))$ is invertible in $\mathbb{R}^{n \times n}$ for any $x \in U$. Therefore the proof of the proposition is complete.  \qed

\begin{coro}
For any (pseudo-) Riemannian metric $g_{ij}[0]$ on $U$, (\ref{gij-q1}) provides a noncommutative metric $g$ on $U$, which is referred as a quantum fluctuation of $g_{ij}[0]$.
\end{coro}
\begin{rmk}
The nonsymmetric metric tensor $g_{ij}[0]$ has its geometric physical origin and the skew-symmetric part represents the intrinsic spin, c.f. \cite{M, KMM, H1, H2} and references therein. The classical gravitational theory in terms of nonsymmetric metric tensor also relates to Connes's noncommutative geometry where spectral triples play roles \cite{KMM}.
\end{rmk}

The noncommutative metric $g$ and its inverse $g^{-1}$ can be written as,
\beQ
\begin{aligned}
    g = & \tilde{E}^i \otimes g_{ij} * E^j = \tilde{E}^i * g_{ij} \otimes E^j,\\
    g^{-1} = & \tilde{E} _i \otimes g^{ij} * E_j = \tilde{E}_i * g^{ij} \otimes E_j.
\end{aligned}
\eeQ
This allows us to define covariant derivatives of $g$ and $g^{-1}$ by
\begin{align*}
\nabla_k g = & \tilde{E}^i \otimes \nabla _k (g_{ij} * E^j) = \tilde{\nabla} _k(\tilde{E}^i * g_{ij}) \otimes E^j,\\
\nabla_k g^{-1} = &\tilde{E} _i \otimes \nabla_k (g^{ij} * E_j) = \tilde{\nabla} _k(\tilde{E}_i * g ^{ij}) \otimes E_j
\end{align*}
and denote
\begin{align*}
\nabla_k g = & \tilde{E}^i \otimes \nabla _k g_{ij} * E^j = \tilde{E}^i * \nabla _k g_{ij} \otimes E^j,\\
\nabla_k g^{-1} = &\tilde{E} _i \otimes \nabla_k g^{ij} * E_j = \tilde{E}_i * \nabla _k g ^{ij} \otimes E_j.
\end{align*}

\begin{prop}\label{nabla-g}
Let $\nabla$ be a noncommutative connection which is compatible with the noncommutative metric $g$. Then
$$ \nabla_k g = \nabla_k g^{-1} =0,$$
i.e.
$$\nabla_k g_{ij} = \nabla_k g^{ij} =0.$$
\end{prop}
\pf It is straightforward that
\begin{align*}
    \nabla_k g &=\nabla_k (\tilde{E}^i \otimes g_{ij} * E^j)\\
    &= \tilde{\nabla}_k\tilde{E}^i \otimes g_{ij} * E^j + \tilde{E}^i \otimes \partial_k g_{ij} * E^j + \tilde{E}^i \otimes g_{ij} * \nabla_kE^j\\
    &= -\tilde{E}^l * \Gamma_{kl}^i \otimes g_{ij} * E^j + \tilde{E}^i \otimes \partial_kg_{ij} * E^j -\tilde{E}^i \otimes g_{ij} * \tilde{\Gamma}_{kl}^j * E^l\\
    &= \tilde{E}^i \otimes (\partial_kg_{ij} - \Gamma_{kij} - \tilde{\Gamma}_{kji}) * E^j =0.
\end{align*}
On the other hand, a direct computation yields
\begin{align*}
    0 &= \Big[\partial_k(g^{il} * g_{lr})\Big]*g^{rj}  \\
      &= \Big[\partial_k{g^{il}} * g_{lr} + g^{il} * \partial_k{g_{lr}} \Big]* g^{rj}\\
      &= \partial_k{g^{il}} * g_{lr} * g^{rj} + g^{il} * (\Gamma_{klr} + \tilde{\Gamma}_{krl}) * g^{rj}\\
      &= \partial_k{g^{ij}} + g^{il} * \Gamma^s_{kl} * g_{sr} * g^{rj} + g^{il} * g_{ls} * \tilde{\Gamma}^s_{kr} * g^{rj}\\
      &= \partial_k{g^{ij}} + g^{il} * \Gamma^j_{kl} + \tilde{\Gamma}^i_{kl} * g^{lj}.
\end{align*}
Therefore,
\begin{align*}
    \nabla_k{g^{-1}} =& \nabla_k (\tilde{E} _i \otimes g^{ij} * E_j)\\
=& \tilde{\nabla}_k\tilde{E}_i \otimes g^{ij} * E_j + \tilde{E}_i \otimes \partial_k g^{ij} * E_j + \tilde{E}_i \otimes g^{ij} * \nabla_k E_j\\
=& \tilde{E}_i \otimes \Big(\partial_k{g^{ij}} + g^{il} * \Gamma^j_{kl} + \tilde{\Gamma}^i_{kl} * g^{lj}\Big) * E_j=0.
\end{align*}\qed

For left and right tangent vectors
\begin{align*}
V=v^i * E_i, \quad W=w^j * E_j, \quad \tilde{V}=\tilde{E} _i * \tilde{v} ^i, \quad \tilde{W}=\tilde{E} _j * \tilde{w} ^j,
\end{align*}
where $v^i$, $w^j$, $\tilde{v}^i$, $\tilde{w}^j  \in \mathcal{A}_U$, the noncommutative Lie brackets are defined as
\begin{align*}
[V, W]f = & v^i * E_i \big(w^j * E_j(f) \big) -  w^j * E_j \big(v^i * E_i(f)\big)\\
      = & \big( v^i * E_i (w^j) -  w^i * E_i (v^j)\big) * E_j (f) +[v^i, w^j] * E_i E_j (f)\\
      = &-[W,V]f,\\
[\tilde{V},\tilde{W}]f=& \tilde{E} _i * \tilde{v}^i  \big(\tilde{E}_j(f) * \tilde{w} ^j\big)
-  \tilde{E}_j * \tilde{w} ^j  \big(\tilde{E}_i (f) *\tilde{v}^i    \big)\\
      = & \tilde{E} _j (f)* \big( \tilde{E} _i (\tilde{w} ^j) * \tilde{v} ^i-  \tilde{E} _i (\tilde{v} ^j)* \tilde{w} ^i\big)
      - \tilde{E} _i \tilde{E} _j (f) * [\tilde{v}^i, \tilde{w}^j]\\
      = &-[\tilde{W},\tilde{V}]f
\end{align*}
for any $f \in \mathcal{A}_U$. Analogous to the classical (pseudo-) Riemannian geometry, the noncommutative left and right curvature operators for left and right tangent vectors can be formally defined as
\begin{align*}
\mathcal{R}_{VW}=&[\nabla_V, \nabla_W] - \nabla_{[V,W]},\\
\mathcal{\tilde{R}}_{\tilde{V}\tilde{W}}=&[\tilde{\nabla}_{\tilde{V}}, \tilde{\nabla}_{\tilde{W}}] - \tilde{\nabla} _{[\tilde{V},\tilde{W}]}.
\end{align*}

It is shown that $\mathcal{R}_{E_iE_j}$, $\tilde{\mathcal{R}}_{\tilde{E}_i\tilde{E}_j}$ are left and right $\mathcal{A}_U$-module endomorphisms over left and right tangent bundles respectively \cite{CTZZ}. But $\mathcal{R}_{VW}$, $\tilde{\mathcal{R}}_{\tilde{V} \tilde{W}}$ do not make sense unless \begin{align*}
[v^i, w^j]=[\tilde{v}^i, \tilde{w}^j]=0.
\end{align*}
Thus, if $V = E_i$ (resp. $\tilde{V} = \tilde{E}_i$) or $W = E_j$ (resp. $\tilde{W} = \tilde{E}_j$), then $\mathcal{R}_{VW}E_k \in \mathcal{T}_U$ (resp. $\tilde{\mathcal{R}}_{\tilde{V}\tilde{W}}\tilde{E}_k \in \tilde{\mathcal{T}}_U$) is well-defined. This suggests to define the covariant derivatives of noncommutative curvatures by adopting the idea of classical (pseudo-) Riemannian geometry. We only consider the case of left curvatures.
\begin{defn} The covariant derivatives of noncommutative curvature operators are defined as follows.
\begin{align*}
    (\nabla_k\mathcal{R})_{E_iE_j}E_p = & \nabla_k(\mathcal{R}_{E_iE_j}E_p) - \mathcal{R}_{(\nabla_{E_k}E_i)E_j}E_p\\
                                        & - \mathcal{R}_{E_i(\nabla_{E_k}E_j)}E_p - \mathcal{R}_{E_iE_j}(\nabla_{E_k}E_p).
\end{align*}
\end{defn}

\begin{defn}
The covariant derivatives of noncommutative curvature tensors, noncommutative Ricci curvatures and noncommutative scalar curvature are defined as follows.
\begin{align*}
    \nabla_s R_{lkij} = & g\big( (\nabla_s \mathcal{R})_{E_iE_j}E_k, \tilde{E}_l \big),\\
    \nabla _s R_j^p = &g^{pk} * g\big( (\nabla_s \mathcal{R})_{E_iE_j}E_k, \tilde{E}_l \big) * g^{li}= g^{pk} * \nabla_s R_{lkij} * g^{li},\\
    \nabla _s \Theta_i^p = & g^{jk} * g\big( (\nabla_s \mathcal{R})_{E_iE_j}E_k, \tilde{E}_l \big) * g^{lp} = g^{jk} * \nabla_s R_{lkij} * g^{lp},\\
    \nabla _s R =& g^{jk} * g\big( (\nabla_s \mathcal{R})_{E_iE_j}E_k, \tilde{E}_l \big) * g^{li}.
\end{align*}
\end{defn}

\begin{rmk}
If $V = E_i$ or $W = E_j$ , then the operator
\begin{align*}
\mathcal{R}_{VW} : \mathcal{T}_U \to \mathcal{T}_U
\end{align*}
is well-defined but generally not $\mathcal{A}_U$-linear. Indeed, let $V = E_i$, $W = a^j * E_j$, we obtain
\begin{align*}
    \mathcal{R}_{VW}(f * E_k) = &\nabla_i\Big(a^j * \nabla_j(f * E_k)\Big) - a^j * \nabla_j\Big(\nabla_i(f * E_k)\Big)\\
    &- \nabla_{[E_i, a^j * E_j]}(f * E_k)\\
    = &(\partial_ia^j) * \nabla_j(f * E_k) + a^j * \nabla_i\Big(\nabla_j(f * E_k)\Big)\\
    &- a^j * \nabla_j\Big(\nabla_i(f * E_k)\Big) - \nabla_{(\partial_ia^j) * E_j}(f * E_k)\\
    = &a^j * \mathcal{R}_{E_iE_j}(f * E_k) + (\partial_ia^j) * \nabla_j(f * E_k)\\
    &- (\partial_ia^j) * \nabla_j(f * E_k)\\
    = &a^j * f * \mathcal{R}_{E_iE_j}E_k.
\end{align*}
The same computation yields that
\begin{align*}
    \mathcal{R}_{VW}E_k = a^j * \mathcal{R}_{E_iE_j}E_k.
\end{align*}
Hence $\mathcal{R}_{VW}(f * E_k)$ and $f * \mathcal{R}_{VW}E_k$ are not equal unless
\begin{align*}
a^j * f = f * a^j.
\end{align*}
The above computation also yields
\begin{align*}
    \mathcal{R}_{E_i(a^j * E_j)} = a^j * \mathcal{R}_{E_iE_j}.
\end{align*}
Similarly, we have
\begin{align*}
    \mathcal{R}_{(a^i * E_i)E_j} = a^i * \mathcal{R}_{E_iE_j}.
\end{align*}

As an operator, $(\nabla_k\mathcal{R})_{E_iE_j}$ dose not give rise to a left $\mathcal{A}_U$-module endomorphism over left tangent bundle. This is because
\begin{align*}
    (\nabla_k\mathcal{R})_{E_iE_j}(f * E_p) = &\nabla_k\Big(\mathcal{R}_{E_iE_j}(f * E_p)\Big) - \mathcal{R}_{(\nabla_kE_i)E_j}(f * E_p)\\
    &- \mathcal{R}_{E_i(\nabla_kE_j)}(f * E_p) - \mathcal{R}_{E_iE_j}\Big(\nabla_k(f * E_p)\Big)\\
    = &\nabla_k(f * \mathcal{R}_{E_iE_j}E_p) - \mathcal{R}_{(\nabla_kE_i)E_j}(f * E_p)\\
    &- \mathcal{R}_{E_i(\nabla_kE_j)}(f * E_p)\\
    &- \mathcal{R}_{E_iE_j}\Big((\partial_kf) * E_p + f * \nabla_kE_p\Big)\\
    = &(\partial_kf) * \mathcal{R}_{E_iE_j}E_p + f * \nabla_k(\mathcal{R}_{E_iE_j}E_p)\\
    &- \mathcal{R}_{(\nabla_kE_i)E_j}(f * E_p) - \mathcal{R}_{E_i(\nabla_kE_j)}(f * E_p)\\
    &- (\partial_kf) * \mathcal{R}_{E_iE_j}E_p - f * \mathcal{R}_{E_iE_j}(\nabla_kE_p)\\
    = &f * \nabla_k(\mathcal{R}_{E_iE_j}E_p) - \mathcal{R}_{(\nabla_kE_i)E_j}(f * E_p)\\
    &- \mathcal{R}_{E_i(\nabla_kE_j)}(f * E_p) - f * \mathcal{R}_{E_iE_j}(\nabla_kE_p)\\
    \not \equiv &f * (\nabla_k\mathcal{R})_{E_iE_j}E_p
\end{align*}
as $\mathcal{R}_{(\nabla_kE_i)E_j}$ and $\mathcal{R}_{E_i(\nabla_kE_j)}$ are not $\mathcal{A}_U$-module endomorphisms in general.
\end{rmk}

\begin{thm}\label{bianchi}
The first (algebraic) Bianchi identity
\begin{align*}
    \mathcal{R}_{E_iE_j}E_k + \mathcal{R}_{E_jE_k}E_i + \mathcal{R}_{E_kE_i}E_j = 0
\end{align*}
and the second (differential) Bianchi identity
\begin{align*}
    (\nabla_i\mathcal{R})_{E_jE_k}E_p + (\nabla_j\mathcal{R})_{E_kE_i}E_p + (\nabla_k\mathcal{R})_{E_iE_j}E_p = 0
\end{align*}
hold for $1 \leq i, j, k, p \leq n$.
\end{thm}
\pf Since the connection is torsion free, we have
\begin{align*}
   \mathcal{R } _{E_iE_j} E_k + & \mathcal{R}_{E_jE_k}E_i + \mathcal{R}_{E_kE_i}E_j \\
     = & \nabla_i\nabla_jE_k - \nabla_j\nabla_iE_k + \nabla_j\nabla_kE_i - \nabla_k\nabla_jE_i \\
       & + \nabla_k\nabla_iE_j - \nabla_i\nabla_kE_j\\
     = & \nabla_i(\nabla_jE_k - \nabla_kE_j)+\nabla_j(\nabla_kE_i - \nabla_iE_k)\\
       & +\nabla_k(\nabla_iE_j - \nabla_jE_i)\\
     = & 0.
\end{align*}
Thus the first Bianchi identity holds. As
\begin{align*}
    (\nabla_i\mathcal{R})_{E_jE_k}E_p = &\nabla_i\nabla_j\nabla_kE_p - \nabla_i\nabla_k\nabla_jE_p - \mathcal{R}_{(\nabla_iE_j)E_k}E_p\\
    &- \mathcal{R}_{E_j(\nabla_iE_k)}E_p - \nabla_j\nabla_k\nabla_iE_p + \nabla_k\nabla_j\nabla_iE_p,\\
    (\nabla_j\mathcal{R})_{E_kE_i}E_p = &\nabla_j\nabla_k\nabla_iE_p - \nabla_j\nabla_i\nabla_kE_p - \mathcal{R}_{(\nabla_jE_k)E_i}E_p\\
    &- \mathcal{R}_{E_k(\nabla_jE_i)}E_p - \nabla_k\nabla_i\nabla_jE_p + \nabla_i\nabla_k\nabla_jE_p,\\
    (\nabla_k\mathcal{R})_{E_iE_j}E_p = &\nabla_k\nabla_i\nabla_jE_p - \nabla_k\nabla_j\nabla_iE_p - \mathcal{R}_{(\nabla_kE_i)E_j}E_p\\
    &- \mathcal{R}_{E_i(\nabla_kE_j)}E_p - \nabla_i\nabla_j\nabla_kE_p + \nabla_j\nabla_i\nabla_kE_p,
\end{align*}
we obtain
\begin{align*}
    &(\nabla_i\mathcal{R})_{E_jE_k}E_p + (\nabla_j\mathcal{R})_{E_kE_i}E_p + (\nabla_k\mathcal{R})_{E_iE_j}E_p\\
    & = -\mathcal{R}_{(\nabla_iE_j)E_k}E_p - \mathcal{R}_{E_j(\nabla_iE_k)}E_p - \mathcal{R}_{(\nabla_jE_k)E_i}E_p\\
    & \quad \,- \mathcal{R}_{E_k(\nabla_jE_i)}E_p - \mathcal{R}_{(\nabla_kE_i)E_j}E_p - \mathcal{R}_{E_i(\nabla_kE_j)}E_p.
\end{align*}
The torsion free condition implies
\begin{align*}
    \mathcal{R}_{(\nabla_iE_j)E_k}E_p + \mathcal{R}_{E_k(\nabla_jE_i)}E_p = 0,\\
    \mathcal{R}_{E_j(\nabla_iE_k)}E_p + \mathcal{R}_{(\nabla_kE_i)E_j}E_p = 0,\\
    \mathcal{R}_{(\nabla_jE_k)E_i}E_p + \mathcal{R}_{E_i(\nabla_kE_j)}E_p = 0.
\end{align*}
Therefore
\begin{align*}
    (\nabla_i\mathcal{R})_{E_jE_k}E_p + (\nabla_j\mathcal{R})_{E_kE_i}E_p + (\nabla_k\mathcal{R})_{E_iE_j}E_p = 0.
\end{align*}
Thus the second Bianchi identity holds. \qed

\begin{rmk}
Bianchi identities also hold for noncommutative right curvature tensors.
\end{rmk}

\begin{prop}
The second Bianchi identity gives that
$$
\nabla_i R_j^i + \nabla_i \Theta_j^i - \delta_j^i   \nabla_i  R=0.
$$
\end{prop}
\pf By the second Bianchi identity, we have
$$
\nabla _i R_{qpjk} +\nabla _j R_{qpki}+\nabla _k R_{qpij}=0.
$$
Multiplying $g^{ip}$ from the left side and $g^{qk}$ from the right side, we obtain
$$
g^{ip}*\nabla _i R_{qpjk}*g^{qk} +g^{ip}*\nabla _j R_{qpki}*g^{qk} + g^{ip}*\nabla _k R_{qpij}*g^{qk}=0.
$$
Taking summation for $i$, $k$, $p$ and $q$, we obtain
$$
-\nabla _i R^i _j +\nabla _j R -\nabla_k \Theta_j ^k=0.
$$
Therefore the proof of the proposition is complete. \qed

\mysection{Equivalence of noncommutative Ricci curvatures}
\ls

In this section, we show that two Ricci curvatures $R^i _j$ and $\Theta ^i _j$ are equivalent under certain conditions. In particular,
they are satisfied if noncommutative metric and chiral coefficients are given by an isometric embedding.

Let noncommutative metric $g_{ij}$, its inverse $g^{ij}$ and chiral coefficients $\Upsilon _{ijk}$ have power series expansions (\ref{gij-q1}), (\ref{gij-q2}) and (\ref{Uijk}).

\begin{lem}\label{lem-inv}
If noncommutative metric $g$ satisfies (\ref{ricci-g}), then
\begin{align}
    g^{ij}[2q] = g^{ji}[2q], \quad g^{ij}[2q + 1] = -g^{ji}[2q + 1]. \label{ricci-gg}
\end{align}
\end{lem}
\pf Since $(g_{ij}[0])$ is symmetric and invertible on $U$, the inverse matrix $(g^{ij}[0])$ is also symmetric, i.e.,
\begin{align*}
    g^{ij}[0] = g^{ji}[0].
\end{align*}
For $u, v \in C^{\infty}(U)$, (\ref{mun}) indicates that
\begin{align}
    \mu_{2q}(u, v) = \mu_{2q}(v, u), \quad \mu_{2q + 1}(u, v) = -\mu_{2q + 1}(v, u).  \label{mun-1}
\end{align}
By (\ref{mun-1}) and the recursive formulas (\ref{r-g}), (\ref{l-g}), we have
\begin{align*}
    g^{ij}[1] &= -g^{ik}[0]g_{kl}[1]g^{lj}[0] - g^{ik}[0]\mu_1 \Big(g_{kl}[0], g^{lj}[0]\Big)\\
    &= g^{jl}[0]g_{lk}[1]g^{ki}[0] + \mu_1 \Big(g^{jl}[0], g_{lk}[0]\Big)g^{ki}[0]\\
    &=-g^{ji}[1].
\end{align*}
Next, let $q \in \mathbb{N}$, using the recursive formula (\ref{r-g}), we have
\begin{align*}
    g^{ij}[2q]
     = &- \sum\limits_{r = 1}^qg^{ik}[0]g_{kl}[2r]g^{lj}[2q - 2r] \\
       &- \sum\limits_{r = 1}^qg^{ik}[0]g_{kl}[2r - 1]g^{lj}[2q - 2r + 1]\\
       &- \sum\limits_{r = 1}^q\sum\limits_{s = 0}^{q - r}g^{ik}[0]\mu_{2r}\Big(g_{kl}[2s], g^{lj}[2q - 2r - 2s]\Big)\\
       &- \sum\limits_{r = 1}^q\sum\limits_{s = 1}^{q - r}g^{ik}[0]\mu_{2r}\Big(g_{kl}[2s - 1], g^{lj}[2q - 2r - 2s + 1]\Big)\\
       &- \sum\limits_{r = 1}^{q}\sum\limits_{s = 0}^{q - r}g^{ik}[0]\mu_{2r - 1}\Big(g_{kl}[2s], g^{lj}[2q - 2r - 2s + 1]\Big)\\
       &- \sum\limits_{r = 1}^q\sum\limits_{s = 0}^{q - r}g^{ik}[0]\mu_{2r - 1}\Big(g_{kl}[2s + 1], g^{lj}[2q - 2r - 2s]\Big).
\end{align*}
Therefore, by induction and (\ref{mun-1}) and recursive formula (\ref{l-g}), we obtain
\begin{align*}
g^{ij}[2q] = &- \sum_{r = 1}^{q}g^{jl}[2q - 2r]g_{lk}[2r]g^{ki}[0] \\
             &- \sum\limits_{r = 1}^q g^{jl}[2q - 2r + 1]g_{lk}[2r - 1]g^{ki}[0]\\
             &- \sum\limits_{r = 1}^q \sum\limits_{s = 0}^{q - r}\mu_{2r}\Big(g^{jl}[2q - 2r - 2s], g_{lk}[2s]\Big)g^{ki}[0]\\
             &- \sum\limits_{r = 1}^q \sum\limits_{s = 1}^{q - r}\mu_{2r}\Big(g^{jl}[2q - 2r - 2s + 1], g_{lk}[2s - 1]\Big)g^{ki}[0]\\
             &- \sum\limits_{r = 1}^q \sum\limits_{s = 0}^{q - r}\mu_{2r - 1}\Big(g^{jl}[2q - 2r - 2s + 1], g_{lk}[2s]\Big)g^{ki}[0]\\
             &- \sum\limits_{r = 1}^q \sum\limits_{s = 0}^{q - r}\mu_{2r - 1}\Big(g^{jl}[2q - 2r - 2s], g_{lk}[2s + 1]\Big)g^{ki}[0]\\
          = & g^{ji}[2q].
\end{align*}
Similarly, we can prove
\begin{align*}
    g^{ij}[2q + 1] = -g^{ji}[2q + 1].
\end{align*}
\qed

\begin{lem}\label{lem-gfg}
If noncommutative metric $g$ satisfies (\ref{ricci-g}), then, for any $f \in C^{\infty}(U)$,
\begin{align*}
    \big(g^{ij} * f * g^{kl} \big)[2q] &= \big(g^{lk} * f * g^{ji} \big)[2q],\\
    \big(g^{ij} * f * g^{kl}\big)[2q + 1] &= -\big(g^{lk} * f * g^{ji}\big)[2q + 1].
\end{align*}
\end{lem}
\pf For any $u, v \in C^{\infty}(U)$, we have
\begin{align*}
    \big(g^{ij} * u * v\big)[2q] = &\sum\limits_{r = 0}^q\sum\limits_{s = 0}^{q - r}\mu_{2r}\Big(g^{ij}[2s], \mu_{2q - 2r - 2s} (u, v)\Big)\\
    &+ \sum\limits_{r = 0}^q\sum\limits_{s = 1}^{q - r}\mu_{2r}\Big(g^{ij}[2s - 1], \mu_{2q - 2r - 2s + 1}(u, v)\Big)\\
    &+ \sum\limits_{r = 1}^q\sum\limits_{s = 0}^{q - r}\mu_{2r - 1}\Big(g^{ij}[2s], \mu_{2q - 2r - 2s + 1}(u, v)\Big)\\
    &+ \sum\limits_{r = 1}^q\sum\limits_{s = 0}^{q - r}\mu_{2r - 1}\Big(g^{ij}[2s + 1], \mu_{2q - 2r - 2s}(u, v)\Big).
\end{align*}
Using Lemma \ref{lem-inv} and (\ref{mun-1}), we obtain
\begin{align*}
\big(g^{ij} * u * v\big)[2q]= &\sum\limits_{r = 0}^q\sum\limits_{s = 0}^{q - r}\mu_{2r}\Big(\mu_{2q - 2r - 2s}(v, u), g^{ji}[2s]\Big)\\
    &+ \sum\limits_{r = 0}^q\sum\limits_{s = 1}^{q - r}\mu_{2r}\Big(\mu_{2q - 2r - 2s + 1}(v, u), g^{ji}[2s - 1]\Big)\\
    &+ \sum\limits_{r = 1}^q\sum\limits_{s = 0}^{q - r}\mu_{2r - 1}\Big(\mu_{2q - 2r - 2s + 1}(v, u), g^{ji}[2s]\Big)\\
    &+ \sum\limits_{r = 1}^q\sum\limits_{s = 0}^{q - r}\mu_{2r - 1}\Big(\mu_{2q - 2r - 2s}(v, u), g^{ji}[2s + 1]\Big)\\
    = &\big(v * u * g^{ji}\big)[2q].
\end{align*}
Similarly, we can show
\begin{align*}
    \big(g^{ij} * u * v\big)[2q + 1] = -\big(v * u * g^{ji}\big)[2q + 1].
\end{align*}
Using them, we obtain
\begin{align*}
    \big(g^{ij} * f * g^{kl}\big)[2q] = &\sum\limits_{r = 0}^q\Big(g^{ij} * f * (g^{kl}[2r])\Big)[2q - 2r]\\
    &+ \sum\limits_{r = 1}^q\Big(g^{ij} * f * (g^{kl}[2r - 1])\Big)[2q - 2r + 1]\\
    = &\sum\limits_{r = 0}^q\Big((g^{kl}[2r]) * f * g^{ji}\Big)[2q - 2r]\\
    &- \sum\limits_{r = 1}^q\Big((g^{kl}[2r - 1]) * f * g^{ji}\Big)[2q - 2r + 1]\\
    = &\sum\limits_{r = 0}^q\Big((g^{lk}[2r]) * f * g^{ji}\Big)[2q - 2r]\\
    &+ \sum\limits_{r = 0}^q\Big((g^{lk}[2r - 1]) * f * g^{ji}\Big)[2q]\\
    = &\big(g^{lk} * f * g^{ji}\big)[2q],
\end{align*}
and, similarly
\begin{align*}
    \big(g^{ij} * f * g^{kl}\big)[2q + 1] = -\big(g^{lk} * f * g^{ji}\big)[2q + 1].
\end{align*}
\qed

\begin{lem}\label{lem-ugv}
If noncommutative metric $g$ satisfies (\ref{ricci-g}), then, for $u, v \in C^{\infty}(U)$,
\begin{align*}
    \big(u * g^{ij} * v\big)[2q] = &\big(v * g^{ji} * u\big)[2q],\\
    \big(u * g^{ij} * v\big)[2q + 1] = &- \big(v * g^{ji} * u\big)[2q + 1].
\end{align*}
\end{lem}
\pf For any $u \in C^{\infty}(U)$, we have
\begin{align*}
    \big(u * g^{ij}\big)[2q] = &\sum\limits_{r = 0}^q\mu_{2r}\Big(u, g^{ij}[2q - 2r]\Big)\\
    &+ \sum\limits_{r = 1}^q\mu_{2r - 1}\Big(u, g^{ij}[2q - 2r + 1]\Big).
\end{align*}
Using Lemma \ref{lem-inv} and (\ref{mun-1}), we obtain
\begin{align*}
    \big(u * g^{ij}\big)[2q] = &\sum\limits_{r = 0}^q\mu_{2r}\Big(u, g^{ij}[2q - 2r]\Big)\\
    &+ \sum\limits_{r = 1}^q\mu_{2r - 1}\Big(u, g^{ij}[2q - 2r + 1]\Big)\\
    = &\sum\limits_{r = 0}^q\mu_{2r}\Big(g^{ji}[2q - 2r], u\Big)\\
    &+ \sum\limits_{r = 1}^q\mu_{2r - 1}\Big(g^{ji}[2q - 2r + 1], u\Big)\\
    = &\big(g^{ji} * u\big)[2q].
\end{align*}
Similarly we have
$$\big(u * g^{ij}\big)[2q + 1] = - \big(g^{ji} * u\big)[2q + 1].$$
Then using Lemma \ref{lem-inv} and (\ref{mun-1}) again, we have
\begin{align*}
    \big(u * g^{ij} * v\big)[2q] = &\sum\limits_{r = 0}^q\mu_{2r}\Big(\big(u * g^{ij}\big)[2q - 2r], v\Big)\\
    &+ \sum\limits_{r = 1}^q\mu_{2r - 1}\Big(\big(u * g^{ij}\big)[2q - 2r + 1], v\Big)\\
    = &\sum\limits_{r = 0}^q\mu_{2r}\Big(v, \big(g^{ji} * u\big)[2q - 2r]\Big)\\
    &+ \sum\limits_{r = 1}^q\mu_{2r - 1}\Big(v, \big(g^{ji} * u\big)[2q - 2r + 1]\Big)\\
    = &\big(v * g^{ji} * u\big)[2q].
\end{align*}
Similarly we have
$$\big(u * g^{ij} * v\big)[2q + 1] = - \big(v * g^{ji} * u\big)[2q + 1].$$\qed

\begin{lem}\label{lem-lrconn}
Let $\nabla$, $\tilde{\nabla}$ be the canonical connection with respect to the noncommutative metric $g$ and chiral coefficients $\Upsilon _{ijk}$ on $U$. If (\ref{ricci-g}), (\ref{ricci-U}) hold, then
\begin{align}\label{con-lr}
\Gamma_{ijk}[2q]  = \tilde{\Gamma}_{ijk}[2q], \quad & \Gamma_{ijk}[2q + 1]  = -\tilde{\Gamma}_{ijk}[2q + 1].
\end{align}
\end{lem}
\pf By chirality and (\ref{ricci-U}), we have
\begin{align*}
    \Gamma_{ijk}[2q] - \tilde{\Gamma}_{ijk}[2q] = &\Upsilon_{ijk}[2q]= 0.
\end{align*}
By (\ref{Gamma-1}), (\ref{Gamma-2}) and (\ref{ricci-g}), we have
\begin{align*}
    \Gamma_{ijk}[2q + 1] = \frac{1}{2}\Upsilon_{ijk}[2q + 1]= - \tilde{\Gamma}_{ijk}[2q + 1].
\end{align*}
\qed

\begin{prop}\label{prop-riem}
Let $M$ be an $n$-dimensional smooth manifold and $U \subset M$ a coordinate chart. Let $\nabla $, $\tilde{\nabla}$ be the canonical connection with respect to noncommutative metric $g$ and chiral coefficients $\Upsilon_{ijk}$ on $U$. Let Riemannian curvatures have the power series expansions (\ref{riem}). If (\ref{ricci-g}), (\ref{ricci-U}) hold, then
\begin{align}
    R_{lkij}[2q] = - R_{klij}[2q], \quad R_{lkij}[2q + 1] = R_{klij}[2q + 1]. \label{1-riem}
\end{align}
\end{prop}
\pf In terms of connection coefficients, the Riemannian curvatures are
\begin{align*}
    R_{lkij} = &g\big((\nabla_i\nabla_j - \nabla_j\nabla_i)E_k, \tilde{E}_l\big)\\
    = &\partial_i g\big(\nabla_jE_k, \tilde{E}_l\big) - g\big(\nabla_jE_k, \tilde{\nabla}_i\tilde{E}_l\big)\\
    &- \partial_j g\big(\nabla_iE_k, \tilde{E}_l\big) + g\big(\nabla_iE_k, \tilde{\nabla}_j\tilde{E}_l\big)\\
    = &\partial_i\Gamma_{jkl} - \partial_j\Gamma_{ikl} + \Gamma^r_{ik} * g_{rs} * \tilde{\Gamma}^s_{jl} - \Gamma^r_{jk} * g_{rs} * \tilde{\Gamma}^s_{il}\\
    = &\partial_i\Gamma_{jkl} - \partial_j\Gamma_{ikl} + \Gamma_{iks} * g^{sr} * \tilde{\Gamma}_{jlr} - \Gamma_{jks} * g^{sr} * \tilde{\Gamma}_{ilr}.
\end{align*}
Since
\begin{align*}
    \partial_i\Gamma_{jkl} - \partial_j\Gamma_{ikl} = &\partial_ig\big(\nabla_jE_k, \tilde{E}_l\big) - \partial_jg\big(\nabla_iE_k, \tilde{E}_l\big)\\
    = &\partial_i\big(\partial_jg(E_k, \tilde{E}_l) - g(E_k, \tilde{\nabla}_j\tilde{E}_l)\big)\\
    &- \partial_j\big(\partial_ig(E_k, \tilde{E}_l) - g(E_k, \tilde{\nabla}_i\tilde{E}_l)\big)\\
    = &\partial_i\partial_jg_{kl} - \partial_i\tilde{\Gamma}_{jlk} - \partial_j\partial_ig_{kl} + \partial_j\tilde{\Gamma}_{ilk}\\
    = &\partial_j\tilde{\Gamma}_{ilk} - \partial_i\tilde{\Gamma}_{jlk},
\end{align*}
we obtain
\begin{align*}
    R_{lkij} = \partial_j\tilde{\Gamma}_{ilk} - \partial_i\tilde{\Gamma}_{jlk} + \Gamma_{iks} * g^{sr} * \tilde{\Gamma}_{jlr} - \Gamma_{jks} * g^{sr} * \tilde{\Gamma}_{ilr}.
\end{align*}
Using Lemma \ref{lem-lrconn}, we have
\begin{align*}
    \big(\partial_i\Gamma_{jkl} - \partial_j\Gamma_{ikl}\big)[2q] = &- \big(\partial_j\tilde{\Gamma}_{ikl} - \partial_i\tilde{\Gamma}_{jkl}\big)[2q]\\
    \big(\partial_i\Gamma_{jkl} - \partial_j\Gamma_{ikl}\big)[2q + 1] = &\big(\partial_j\tilde{\Gamma}_{ikl} - \partial_i\tilde{\Gamma}_{jkl}\big)[2q + 1].
\end{align*}
Moreover, using Lemma \ref{lem-ugv} and Lemma \ref{lem-lrconn}, we obtain
\begin{align*}
    \mathrel{\phantom{=}}&\big(\Gamma_{iks} * g^{sr} * \tilde{\Gamma}_{jlr}\big)[2q]\\
    = &\sum\limits_{\alpha = 0}^q\sum\limits_{\beta = 0}^{q - \alpha}\big(\Gamma_{iks}[2\alpha] * g^{sr} * \tilde{\Gamma}_{jlr}[2\beta]\big)[2q - 2\alpha - 2\beta]\\
    &+ \sum\limits_{\alpha = 0}^q\sum\limits_{\beta = 1}^{q - \alpha}\big(\Gamma_{iks}[2\alpha] * g^{sr} * \tilde{\Gamma}_{jlr}[2\beta - 1]\big)[2q - 2\alpha - 2\beta + 1]\\
    &+ \sum\limits_{\alpha = 1}^q\sum\limits_{\beta = 0}^{q - \alpha}\big(\Gamma_{iks}[2\alpha - 1] * g^{sr} * \tilde{\Gamma}_{jlr}[2\beta]\big)[2q - 2\alpha - 2\beta + 1]\\
    &+ \sum\limits_{\alpha = 1}^q\sum\limits_{\beta = 0}^{q - \alpha}\big(\Gamma_{iks}[2\alpha - 1] * g^{sr} * \tilde{\Gamma}_{jlr}[2\beta + 1]\big)[2q - 2\alpha - 2\beta]\\
    = &\sum\limits_{\alpha = 0}^q\sum\limits_{\beta = 0}^{q - \alpha}\big(\Gamma_{jlr}[2\beta] * g^{rs} * \tilde{\Gamma}_{iks}[2\alpha]\big)[2q - 2\alpha - 2\beta]\\
    &+ \sum\limits_{\alpha = 0}^q\sum\limits_{\beta = 1}^{q - \alpha}\big(\Gamma_{jlr}[2\beta - 1] * g^{rs} * \tilde{\Gamma}_{iks}[2\alpha]\big)[2q - 2\alpha - 2\beta + 1]\\
    &+ \sum\limits_{\alpha = 1}^q\sum\limits_{\beta = 0}^{q - \alpha}\big(\Gamma_{jlr}[2\beta] * g^{rs} * \tilde{\Gamma}_{iks}[2\alpha - 1]\big)[2q - 2\alpha - 2\beta + 1]\\
    &+ \sum\limits_{\alpha = 1}^q\sum\limits_{\beta = 0}^{q - \alpha}\big(\Gamma_{jlr}[2\beta + 1] * g^{rs} * \tilde{\Gamma}_{iks}[2\alpha - 1]\big)[2q - 2\alpha - 2\beta]\\
    = &\big(\Gamma_{jls} * g^{sr} * \tilde{\Gamma}_{ikr}\big)[2q].
\end{align*}
Similarly, we have
\begin{align*}
    \big(\Gamma_{iks} * g^{sr} * \tilde{\Gamma}_{jlr}\big)[2q + 1] = &- \big(\Gamma_{jls} * g^{sr} * \tilde{\Gamma}_{ikr}\big)[2q + 1],\\
    \big(\Gamma_{jks} * g^{sr} * \tilde{\Gamma}_{ilr}\big)[2q] = &\big(\Gamma_{ils} * g^{sr} * \tilde{\Gamma}_{jkr}\big)[2q],\\
    \big(\Gamma_{jks} * g^{sr} * \tilde{\Gamma}_{ilr}\big)[2q + 1] = &- \big(\Gamma_{ils} * g^{sr} * \tilde{\Gamma}_{jkr}\big)[2q + 1].
\end{align*}
Therefore
\begin{align*}
    R_{lkij}[2q] = &\big(\partial_i\Gamma_{jkl} - \partial_j\Gamma_{ikl} + \Gamma_{iks} * g^{sr} * \tilde{\Gamma}_{jlr} - \Gamma_{jks} * g^{sr} * \tilde{\Gamma}_{ilr}\big)[2q]\\
    = &- \big(\partial_j\tilde{\Gamma}_{ikl} - \partial_i\tilde{\Gamma}_{jkl} + \Gamma_{ils} * g^{sr} * \Gamma_{jkr} - \Gamma_{jls} * g^{sr} * \tilde{\Gamma}_{ikr}\big)[2q]\\
    = &- R_{klij}[2q].
\end{align*}
Similarly, we obtain
$$R_{lkij}[2q + 1] = R_{klij}[2q + 1].$$
\qed

\begin{thm}\label{ricci}
Let $M$ be an $n$-dimensional smooth manifold and $U \subset M$ a coordinate chart. Let $\nabla $, $\tilde{\nabla}$ be the canonical connection with respect to noncommutative metric $g$ and chiral coefficients $\Upsilon_{ijk}$ on $U$. Let two Ricci curvatures have the power series expansions (\ref{lowricci-q1}), (\ref{lowricci-q2}), (\ref{ricci-q1}), (\ref{ricci-q2}). If (\ref{ricci-g}), (\ref{ricci-U}) hold, then
\begin{align}
    R_{ij}[2q] &= \Theta_{ji}[2q], \quad R_{ij}[2q + 1] = - \Theta_{ji}[2q + 1], \label{1-ricci}\\
    R^i_j[2q] &= \Theta^i_j[2q], \quad   \,\,\,R^i_j[2q + 1] = -\Theta^i_j[2q + 1]. \label{2-ricci}
\end{align}
\end{thm}
\pf By definition
\begin{align*}
    R_{ij} = R_{likj} * g^{lk}, \quad \Theta_{ij} = g^{lk} * R_{jkil}.
\end{align*}
Using Lemma \ref{lem-inv}, Proposition \ref{prop-riem} and (\ref{mun-1}), we obtain
\begin{align*}
    R_{ij}[2q] = &\sum\limits_{r = 0}^q\sum\limits_{s = 0}^{q - r}\mu_{2r}\Big(R_{likj}[2s], g^{lk}[2q - 2r - 2s]\Big)\\
    &+ \sum\limits_{r = 0}^q\sum\limits_{s = 1}^{q - r}\mu_{2r}\Big(R_{likj}[2s - 1], g^{lk}[2q - 2r - 2s + 1]\Big)\\
    &+ \sum\limits_{r = 1}^q\sum\limits_{s = 0}^{q - r}\mu_{2r - 1}\Big(R_{likj}[2s], g^{lk}[2q - 2r - 2s + 1]\Big)\\
    &+ \sum\limits_{r = 1}^q\sum\limits_{s = 0}^{q - r}\mu_{2r - 1}\Big(R_{likj}[2s + 1], g^{lk}[2q - 2r - 2s]\Big)\\
    = &- \sum\limits_{r = 0}^q\sum\limits_{s = 0}^{q - r}\mu_{2r}\Big(g^{kl}[2q - 2r - 2s], R_{ilkj}[2s]\Big)\\
    &- \sum\limits_{r = 0}^q\sum\limits_{s = 1}^{q - r}\mu_{2r}\Big(g^{kl}[2q - 2r - 2s + 1], R_{ilkj}[2s - 1]\Big)\\
    &- \sum\limits_{r = 1}^q\sum\limits_{s = 0}^{q - r}\mu_{2r - 1}\Big(g^{kl}[2q - 2r - 2s + 1], R_{ilkj}[2s]\Big)\\
    &- \sum\limits_{r = 1}^q\sum\limits_{s = 0}^{q - r}\mu_{2r - 1}\Big(g^{kl}[2q - 2r - 2s], R_{ilkj}[2s + 1]\Big)\\
    = &- \big(g^{kl} * R_{ilkj}\big)[2q]\\
    = &\big(g^{kl} * R_{iljk}\big)[2q]\\
    = &\Theta_{ji}[2q].
\end{align*}
Similarly, we have
$$R_{ij}[2q + 1] = - \Theta_{ji}[2q + 1].$$
Since
$$R^i_j = g^{ik} * R_{kj}, \quad \Theta^i_j = \Theta_{jk} * g^{ki},$$
we obtain
\begin{align*}
    R^i_j[2q] = &\sum\limits_{r = 0}^q\sum\limits_{s = 0}^{q - r}\mu_{2r}\Big(g^{ik}[2s], R_{kj}[2q - 2r - 2s]\Big)\\
    &+ \sum\limits_{r = 0}^q\sum\limits_{s = 1}^{q - r}\mu_{2r}\Big(g^{ik}[2s - 1], R_{kj}[2q - 2r - 2s + 1]\Big)\\
    &+ \sum\limits_{r = 1}^q\sum\limits_{s = 0}^{q - r}\mu_{2r - 1}\Big(g^{ik}[2s], R_{kj}[2q - 2r - 2s + 1]\Big)\\
    &+ \sum\limits_{r = 1}^q\sum\limits_{s = 0}^{q - r}\mu_{2r - 1}\Big(g^{ik}[2s + 1], R_{kj}[2q - 2r - 2s]\Big)\\
    = &\sum\limits_{r = 0}^q\sum\limits_{s = 0}^{q - r}\mu_{2r}\Big(\Theta_{jk}[2q - 2r - 2s], g^{ki}[2s]\Big)\\
    &+ \sum\limits_{r = 0}^q\sum\limits_{s = 1}^{q - r}\mu_{2r}\Big(\Theta_{jk}[2q - 2r - 2s + 1], g^{ki}[2s - 1]\Big)\\
    &+ \sum\limits_{r = 1}^q\sum\limits_{s = 0}^{q - r}\mu_{2r - 1}\Big(\Theta_{jk}[2q - 2r - 2s + 1], g^{ki}[2s]\Big)\\
    &+ \sum\limits_{r = 1}^q\sum\limits_{s = 1}^{q - r}\mu_{2r - 1}\Big(\Theta_{jk}[2q - 2r - 2s], g^{ki}[2s + 1]\Big)\\
    = &\big(\Theta_{jk} * g^{ki}\big)[2q]\\
    = &\Theta^i_j[2q].
\end{align*}
Similarly, we have
$$R^i_j[2q + 1] = - \Theta^i_j[2q + 1].$$
\qed

\begin{rmk}
We provide two examples on $\mathbb{R}^2$. The metric $g$ does not satisfy (\ref{ricci-g}) in the first example and the chiral coefficients does not satisfy (\ref{ricci-U}) in the second example. In both examples the two Ricci curvatures do not satisfy (\ref{1-ricci}) and (\ref{2-ricci}).

Example 1. Let the metric
\begin{equation*}
    g = \left(
    \begin{array}{cc}
    \sum\limits_{q = 0}^{\infty}\hbar^q & 0 \\
    0 & \sum\limits_{q = 0}^{\infty}\hbar^q
    \end{array}
    \right),
\end{equation*}
the chiral coefficients satisfy
\begin{align*}
    \Upsilon_{111} &= \hbar, &\Upsilon_{121}& = \Upsilon_{211} = - \hbar, &\Upsilon_{221}& = - \hbar,\\
    \Upsilon_{112} &= - \hbar, &\Upsilon_{122}& = \Upsilon_{212} = - \hbar, &\Upsilon_{222}& = \hbar.
\end{align*}

It is straightforward that the connection coefficients are
\begin{align*}
    \Gamma_{111} &= \frac{1}{2}\hbar, &\Gamma_{121}& = \Gamma_{211} = - \frac{1}{2}\hbar, &\Gamma_{221}& = - \frac{1}{2}\hbar,\\
    \Gamma_{112} &= - \frac{1}{2}\hbar, &\Gamma_{122}& = \Gamma_{212} = - \frac{1}{2}\hbar, &\Gamma_{222}& = \frac{1}{2}\hbar,\\
    \tilde{\Gamma}_{111} &= - \frac{1}{2}\hbar, &\tilde{\Gamma}_{121}& = \tilde{\Gamma}_{211} = \frac{1}{2}\hbar, &\tilde{\Gamma}_{221}& = \frac{1}{2}\hbar,\\
    \tilde{\Gamma}_{112} &= \frac{1}{2}\hbar, &\tilde{\Gamma}_{122}& = \tilde{\Gamma}_{212} = \frac{1}{2}\hbar, &\tilde{\Gamma}_{222}& = - \frac{1}{2}\hbar
\end{align*}
and the inverse metric is
\begin{equation*}
    g^{- 1} = \left(
    \begin{array}{cc}
    1 - \hbar & 0 \\
    0 & 1 - \hbar
    \end{array}
    \right).
\end{equation*}
Then the nontrivial components of the Riemannian and Ricci curvatures are
\begin{align*}
    R_{1212} &= - R_{1221} = - R_{2112} = R_{2121} = - \hbar^2(1 - \hbar),\\
    R_{11} &= R_{22} = \Theta_{11} = \Theta_{22} = - \hbar^2(1 - \hbar)^2,\\
    R^1_1 &= R^2_2 = \Theta^1_1 = \Theta^2_2 = - \hbar^2(1 - \hbar)^3.
\end{align*}
Therefore (\ref{1-ricci}) and (\ref{2-ricci}) fail.

Example 2. Let the metric
\begin{equation*}
    g = \left(
    \begin{array}{cc}
    \sum\limits_{q = 0}^{\infty}(-1)^q\hbar^{2q} & 0 \\
    0 & \sum\limits_{q = 0}^{\infty}(-1)^q\hbar^{2q}
    \end{array}
    \right),
\end{equation*}
the chiral coefficients satisfy
\begin{align*}
    \Upsilon_{111} &= \hbar, &\Upsilon_{121}&  = \Upsilon_{211} = - \hbar, &\Upsilon_{221}& = - \hbar,\\
    \Upsilon_{112} &= - \hbar, &\Upsilon_{122}&  = \Upsilon_{212} = - \hbar, &\Upsilon_{222}& = \hbar^2.
\end{align*}

It is straightforward that the connection coefficients are
\begin{align*}
    \Gamma_{111} &= \frac{1}{2}\hbar, &\Gamma_{121}& = \Gamma_{211} = - \frac{1}{2}\hbar, &\Gamma_{221}& = - \frac{1}{2}\hbar,\\
    \Gamma_{112} &= - \frac{1}{2}\hbar, &\Gamma_{122}& = \Gamma_{212} = - \frac{1}{2}\hbar, &\Gamma_{222}& = \frac{1}{2}\hbar^2,\\
    \tilde{\Gamma}_{111} &= - \frac{1}{2}\hbar, &\tilde{\Gamma}_{121}& = \tilde{\Gamma}_{211} = \frac{1}{2}\hbar,  &\tilde{\Gamma}_{221}& = \frac{1}{2}\hbar,\\
    \tilde{\Gamma}_{112} &= \frac{1}{2}\hbar, &\tilde{\Gamma}_{122}& = \tilde{\Gamma}_{212} = \frac{1}{2}\hbar, &\tilde{\Gamma}_{222}& = - \frac{1}{2}\hbar^2
\end{align*}
and the inverse metric is
\begin{equation*}
    g^{- 1} = \left(
    \begin{array}{cc}
    1 + \hbar^2 & 0 \\
    0 & 1 + \hbar^2
    \end{array}
    \right).
\end{equation*}
Then the nontrivial components of the Riemannian and Ricci curvatures are
\begin{align*}
    R_{1212} &= - R_{1221} = - R_{2112} = R_{2121} = - \Big(\frac{3}{4}\hbar^2 + \frac{1}{4}\hbar^3 + \frac{3}{4}\hbar^4 + \frac{1}{4}\hbar^5\Big),\\
    R_{11} &= R_{22} = \Theta_{11} = \Theta_{22} = - \Big(1 + \hbar^2\Big)\Big(\frac{3}{4}\hbar^2 + \frac{1}{4}\hbar^3 + \frac{3}{4}\hbar^4 + \frac{1}{4}\hbar^5\Big),\\
    R^1_1 &= R^2_2 = \Theta^1_1 = \Theta^2_2 = - \Big(1 + \hbar^2\Big)^2\Big(\frac{3}{4}\hbar^2 + \frac{1}{4}\hbar^3 + \frac{3}{4}\hbar^4 + \frac{1}{4}\hbar^5\Big).
\end{align*}
Therefore (\ref{1-ricci}) and (\ref{2-ricci}) fail.
\end{rmk}

Now we provide the quantum fluctuation of a pseudo-Riemannian metric $g[0]$ on $U$ in terms of isometric embedding \cite{CTZZ, WZZ, WZZ1}. Recall that $(U, g[0])$ can always be isometrically embedded into a pseudo-Euclidean space, c.f. \cite{PT}, i.e., there exist a differentiable map
$$
X: U \lrw \mathbb{R}^{p, m - p}
$$
such that
$$
g_{ij}[0]=\sum\limits_{\alpha = 1}^m \eta _{\alpha \alpha} \partial_i X^\alpha \cdot \partial _j X^\alpha,
$$
where $\eta =\mbox{diag} (-1,\cdots,-1, 1,\cdots,1)$ is the flat metrics of $\mathbb{R}^{p, m - p}$. The quantum fluctuation of $g[0]$ is
\begin{align}
g\big(E_i, \tilde{E}_j \big)=\sum\limits_{\alpha = 1}^m \eta _{\alpha\alpha} \partial_i X^\alpha * \partial _j X^\alpha, \label{g-iso}
\end{align}
where $E_i =\tilde{E}_i =\partial_i$. It yields a canonical connection with the connection and chiral coefficients
\begin{align}
\Gamma _{ijk} =&\sum\limits_{\alpha = 1}^m \eta _{\alpha \alpha} \partial_i \partial_j X^\alpha * \partial _k X^\alpha, \label{Gamma-iso1}\\
\tilde{\Gamma} _{ijk} =&\sum\limits_{\alpha = 1}^m \eta _{\alpha \alpha} \partial _k X^\alpha * \partial_i \partial_j X^\alpha, \label{Gamma-iso2}\\
\Upsilon_{ijk} =&\sum\limits_{\alpha = 1}^m \eta _{\alpha \alpha} \big(\partial_i \partial_j X^\alpha * \partial _k X^\alpha -\partial_k X^\alpha * \partial_i \partial_j X^\alpha \big). \label{U-iso}
\end{align}

\begin{coro}\label{ricci-X}
Let $g$ be given by (\ref{g-iso}) for isometric embedding
$$X = (X^1, \cdots, X^m) \in C^{\infty}(U, \mathbb{R}^m),$$ where $X^{\alpha} \in C^{\infty}(U)$, $1 \leq \alpha \leq m$. Let $R_j^i$ and $\Theta_j^i$ be the two Ricci curvatures of the canonical connection induced by $X$. Then (\ref{1-ricci}) and (\ref{2-ricci}) hold.
\end{coro}
\pf We only need to check (\ref{g-iso}), (\ref{U-iso}) satisfy (\ref{ricci-g}) and (\ref{ricci-U}). Indeed, it is a direct consequence of (\ref{mun-1}). \qed

\begin{rmk}
The following noncommutative Einstein field equations were proposed in \cite{CTZZ}
$$
R_j^i + \Theta_j^i - \delta_j^iR =T^i _j.
$$
As it may not capture all information of noncommutative metrics, the second author gave the strong version in \cite{Z}
$$
   R_j^i - \frac{1}{2}\delta_j^iR = T_j^i, \quad    \Theta_j^i - \frac{1}{2}\delta_j^iR = \tilde{T}_j^i.
$$
Theorem \ref{ricci} and Corollary \ref{ricci-X} indicate that only the first one is sufficient and the noncommutative Einstein field equations should be
\begin{align*}
    R_j^i - \frac{1}{2}\delta_j^iR = T_j^i
\end{align*}
if (\ref{ricci-g}), (\ref{ricci-U}) hold, in particular, if noncommutative metrics are given by isometric embedding.
\end{rmk}

\mysection{Spherically symmetric isometric embedding}
\ls

In this section, we show that the quantum fluctuations and their curvatures have closed forms coming from Moyal products of trigonometric functions if (pseudo-) Riemannian metrics are given by certain type of spherically symmetric isometric embedding. This indicates that the quantization of gravity is renormalizable in this case.

\begin{thm}\label{renor}
Let open set
\begin{align*}
U = (0, \infty) \times (0, 2\pi) \times (0, \pi) \times \cdots \times (0, \pi) \subset \mathbb{R}^{n},
\end{align*}
which is equipped with coordinates $(x^1, x^2, \cdots, x^n) = (\rho, \theta_1, \cdots, \theta_{n - 1})$. Let $(U, g[0])$ be a (pseudo-) Riemannian metric given by a spherically symmetric isometric embedding
$$X : U \lrw \mathbb{R}^{p, m - p}$$
with
\begin{align*}
    X^1           &= f^1(\rho),\\
                  & \cdots \cdots\\
    X^{m - n}     &= f^{m - n}(\rho),\\
    X^{m - n + 1} &= f^{m - n + 1}(\rho)\sin\theta_{n - 1} \sin\theta_{n - 2}{\cdots} \sin\theta_2 \sin\theta_1,\\
    X^{m - n + 2} &= f^{m - n + 2}(\rho)\sin\theta_{n - 1} \sin\theta_{n - 2}{\cdots} \sin\theta_2 \cos\theta_1,\\
                  & \cdots \cdots\\
    X^{m - 2}     &= f^{m - 2}(\rho)\sin\theta_{n - 1} \sin\theta_{n - 2} \cos\theta_{n - 3},\\
    X^{m - 1}     &= f^{m - 1}(\rho)\sin\theta_{n - 1} \cos\theta_{n - 2},\\
    X^m           &= f^m(\rho)\cos \theta_{n - 1},
\end{align*}
where $f^1(\rho)$, $\cdots$, $f^m(\rho)$ are smooth functions of $\rho$, $m - n + 1 > p$ and
\begin{align*}
f^{m - n + 1}(\rho)= f^{m - n + 2}(\rho)=f(\rho).
\end{align*}
Fix some $l \in [3, n]$, define the Moyal product in terms of skew-symmetric matrix $(\theta^{ij})$ with nonzero elements
$$\theta^{2 l} = - \theta^{l 2} = \lambda \neq 0.$$
Then the quantum fluctuation of $g[0]$ and their curvatures have closed forms coming from absolutely convergent power series expansions on $U$.
\end{thm}
\pf Note that only the term
$$\partial_i X^{m - n + 1} * \partial_j X^{m - n + 1} + \partial_iX^{m - n + 2} * \partial_jX^{m - n + 2}$$
cannot be reduced to the usual commutative product in noncommutative metric (\ref{g-iso}). Denote
$$g_{ij}^{(\alpha, \beta)} = \partial_i X^\alpha * \partial_j X^\beta.$$
Denote $a_0=m-n+1$ for short. Using the formulas provided in the appendix, we obtain, for $2 < i<j \leq n, 2 < k \leq n$ and $i, j, k \neq l$,
\begin{align*}
g_{11}^{(a_0, a_0)} + g_{11}^{(a_0+1, a_0+1)}
    = &\big(f'\big)^2 \sin^2 \theta_{n - 1} \cdots \sin^2 \theta_l \sin^2 \theta_{l - 2} \cdots \sin^2 \theta_2\\
      & \Big(\sin^2 \theta_{l - 1} \cosh^2(\lambda \hbar) - \cos^2 \theta_{l - 1} \sinh^2(\lambda \hbar)\Big),\\
g_{12}^{(a_0, a_0)} + g_{12}^{(a_0+1, a_0+1)} = & - \Big(g_{21}^{(a_0, a_0)} + g_{21}^{(a_0+1, a_0+1)}\Big)\\
    = & 2 f f'\sin^2 \theta_{n - 1} \cdots \sin^2 \theta_l \sin^2 \theta_{l - 2} \cdots \sin^2 \theta_2\\
      & \sin \theta_{l - 1} \cos \theta_{l - 1} \cosh(\lambda \hbar)\sinh(\lambda \hbar),\\
g_{1l}^{(a_0, a_0)} + g_{1l}^{(a_0+1, a_0+1)} = & g_{l1}^{(a_0, a_0)} + g_{l1}^{(a_0+1, a_0+1)}\\
    = & f f' \sin^2 \theta_{n - 1} \cdots \sin^2 \theta_l \sin^2 \theta_{l - 2} \cdots \sin^2 \theta_2\\
      & \sin \theta_{l - 1} \cos \theta_{l - 1} \Big(1 + 2\sinh^2(\lambda \hbar)\Big),\\
g_{22}^{(a_0, a_0)} + g_{22}^{(a_0+1, a_0+1)}
    = & f^2 \sin^2 \theta_{n - 1} \cdots \sin^2 \theta_l \sin^2 \theta_{l - 2} \cdots \sin^2 \theta_2\\
      & \Big(\sin^2 \theta_{l - 1} \cosh^2(\lambda \hbar) - \cos^2 \theta_{l - 1} \sinh^2(\lambda \hbar)\Big),\\
g_{2l}^{(a_0, a_0)} + g_{2l}^{(a_0+1, a_0+1)} = & - \Big(g_{l2}^{(a_0, a_0)} + g_{l2}^{(a_0+1, a_0+1)}\Big)\\
    = & f^2 \sin^2 \theta_{n - 1} \cdots \sin^2 \theta_l \sin^2 \theta_{l - 2} \cdots \sin^2 \theta_2\\
      & \Big(\sin^2 \theta_{l - 1} - \cos^2 \theta_{l - 1} \Big) \cosh(\lambda \hbar) \sinh(\lambda \hbar),\\
g_{ll}^{(a_0, a_0)} + g_{ll}^{(a_0+1, a_0+1)}
    = & f ^2 \sin^2 \theta_{n - 1} \cdots \sin^2 \theta_l \sin^2 \theta_{l - 2} \cdots \sin^2 \theta_2\\
      & \Big(\cos^2 \theta_{l - 1} \cosh^2(\lambda \hbar) - \sin^2 \theta_{l - 1} \sinh^2(\lambda \hbar)\Big),\\
g_{1k}^{(a_0, a_0)} + g_{1k}^{(a_0+1, a_0+1)} = & g_{k1}^{(a_0, a_0)} + g_{k1}^{(a_0+1, a_0+1)}\\
    = &f f' \sin^2 \theta_{n - 1} \cdots \sin^2 \theta_k \sin \theta_{k - 1} \cos \theta_{k - 1} \sin^2 \theta_{k - 2}\\
      & \cdots\Big(\sin^2 \theta_{l - 1} \cosh^2(\lambda \hbar) - \cos^2 \theta_{l - 1} \sinh^2(\lambda \hbar)\Big),\\
g_{2k}^{(a_0, a_0)} + g_{2k}^{(a_0+1, a_0+1)} = & - \Big(g_{k2}^{(a_0, a_0)} + g_{k2}^{(a_0+1, a_0+1)}\Big)\\
    = &- 2f^2 \sin^2 \theta_{n - 1} \cdots \sin^2 \theta_k \sin \theta_{k - 1} \cos \theta_{k - 1} \sin^2 \theta_{k - 2}\\
    & \cdots\sin \theta_{l - 1} \cos \theta_{l - 1} \cosh(\lambda \hbar) \sinh(\lambda \hbar),\\
g_{lk}^{(a_0, a_0)} + g_{lk}^{(a_0+1, a_0+1)} = & g_{kl}^{(a_0, a_0)} + g_{kl}^{(a_0+1, a_0+1)}\\
    = & f^2 \sin^2 \theta_{n - 1} \cdots \sin^2 \theta_k \sin \theta_{k - 1} \cos \theta_{k - 1} \sin^2 \theta_{k - 2}\\
      & \cdots\sin \theta_{l - 1} \cos \theta_{l - 1} \Big(1 + 2\sinh^2(\lambda \hbar)\Big),\\
      g_{kk}^{(a_0, a_0)} + g_{kk}^{(a_0+1, a_0+1)}
    = & f^2 \sin^2 \theta_{n - 1} \cdots \sin^2 \theta_k \cos^2 \theta_{k - 1} \sin^2 \theta_{k - 2}\\
      & \cdots\Big(\sin^2 \theta_{l - 1} \cosh^2(\lambda \hbar) - \cos^2 \theta_{l - 1} \sinh^2(\lambda \hbar)\Big),\\
g_{ij}^{(a_0, a_0)} + g_{ij}^{(a_0+1, a_0+1)} = & g_{ji}^{(a_0, a_0)} + g_{ji}^{(a_0+1, a_0+1)}\\
    = & f^2 \sin^2 \theta_{n - 1} \cdots \sin^2 \theta_j \sin \theta_{j - 1} \cos \theta_{j - 1} \sin^2 \theta_{j - 2}\\
      & \cdots\sin^2 \theta_i \sin \theta_{i - 1} \cos \theta_{i - 1} \sin^2 \theta_{i - 2}\\
      & \cdots\Big(\sin^2 \theta_{l - 1} \cosh^2(\lambda \hbar) - \cos^2 \theta_{l - 1} \sinh^2(\lambda \hbar)\Big).
\end{align*}
They indicate that
\begin{align*}
g_{2k}=-g_{k2}, \quad k \neq 2
\end{align*}
but other metric components are symmetric, and the quantum fluctuation $g = (g_{ij})$ of $g[0]$ have closed forms which are smooth functions not depending on $\theta_1$. Therefore the Moyal product relating to $(g_{ij})$ becomes usual commutative product. This means the inverse matrix $(g^{ij})$ coincides with the inverse matrix in the sense of usual commutative product, and its elements do not depend on $\theta _1$ neither. By (\ref{Gamma-iso1}), (\ref{Gamma-iso2}), similar calculation yields that the connection coefficients $\Gamma_{ijk}$, $\tilde{\Gamma}_{ijk}$ also have closed forms which are smooth functions not depending on $\theta_1$.

As all quantities relating to the quantum fluctuation and the connection coefficients do not depend on $\theta _1$, the Moyal product in deriving the curvatures becomes usual commutative product. Therefore the curvatures have closed forms, which depend only on $\rho$, $\theta_2, \cdots, \theta_{n - 1}$ and $\hbar$.
\qed

\mysection{Quasi-connections and curvatures}
\ls

In this section, we study noncommutative differential geometry with respect to star products constructed by Kontsevich on Poisson manifolds \cite{Ko}. These star products are only associative, but are not compatible with the Leibniz rule. This causes that many geometric properties are lost. However, we can still define left and right quasi-connections as well as their curvatures.

Let $M$ be an $n$-dimensional differentiable manifold and $U \subset M$ a coordinate chart equipped with coordinates $\{ x^1, \cdots, x^n \}$. Recall that a formal deformation of the $\mathbb{R}$-algebra $C^{\infty}(U)$ is an associative $\mathbb{R}[[\hbar]]$-bilinear product
\begin{align*}
\star : \mathcal{A}_U \times \mathcal{A}_U \longrightarrow \mathcal{A}_U
\end{align*}
such that, for $u,\,\, v \in C^{\infty}(U)$,
\begin{align*}
1 \star u = u \star 1 = u,\quad u \star v = uv + \sum\limits_{q = 1}^{\infty}B_q(u, v)\hbar^q,
\end{align*}
where
\begin{align*}
B_q : C^{\infty}(U) \times C^{\infty}(U) \longrightarrow C^{\infty}(U)
\end{align*}
are $\mathbb{R}$-bilinear maps. We call $\star$ a deformation product on $C^{\infty}(U)$ and $(\mathcal{A}_U, \star)$ a deformation algebra of $C^{\infty}(U)$. The associativity of $\star$ implies that
\begin{align*}
\{u, v\}_{\star}=\frac{1}{2} \Big(B_1(u, v) - B_1(v, u)\Big)
\end{align*}
is a Poisson bracket, c.f. \cite{Ko, LPV}. Furthermore, we call $\star$ a star product if all $B_q$ are bi-differential operators. Given a Poisson bracket $\{, \}$ on $C^{\infty}(U)$, a star product $\star$ is called a deformation quantization of $\big(U, \{, \}\big)$ if
\begin{align*}
\{, \} = \{, \}_{\star},
\end{align*}
c.f. \cite{BFFLS2}. The existence of deformation quantization on Poisson manifolds was proved by Kontsevich \cite{Ko}.

As before, we denote
\begin{align*}
E_i = \tilde{E}_i = \partial_i, \quad 1 \leq i \leq n.
\end{align*}
The noncommutative tangent bundles and metrics with respect to a star products can be defined as follows.
\begin{defn}
The noncommutative left ({\rm{resp.}} right) $\star$-tangent bundle $^{\star}\mathcal{T}_U$ ({\rm{resp.}} $^{\star}\tilde{\mathcal{T}}_U$) on $U$ is the free left ({\rm{resp.}} right) $(\mathcal{A}_U, \star)$-module with basis $\{E_1, \cdots, E_n\}$ ({\rm{resp.}} $\{\tilde{E}_1, \cdots, \tilde{E}_n\}$), i.e.,
\begin{align*}
^{\star}\mathcal{T}_U & = \Big\{a^i \star E_i \,\Big | \,a^i \in \mathcal{A}_U, \, a^i \star E_i = 0 \Longleftrightarrow a^i = 0 \Big\},\\
^{\star}\tilde{\mathcal{T}}_U & = \Big\{\tilde{E}_i \star \tilde{a}^i \,\Big | \, \tilde{a}^i \in \mathcal{A}_U, \, \tilde{E}_i \star \tilde{a}^i = 0 \Longleftrightarrow \tilde{a}^i = 0 \Big\}.
\end{align*}
An element of $^{\star}\mathcal{T}_U$ ({\rm{resp.}} $^{\star}\tilde{\mathcal{T}}_U$) is called a left ({\rm{resp.}} right) $\star$-vector field.
\end{defn}

\begin{defn}
A noncommutative $\star$-metric $^{\star}g$ is defined as a homomorphism of two-sided $(\mathcal{A}_U, \star)$-modules
\begin{align*}
^{\star}g : {^{\star}\mathcal{T}_U} \otimes_{\mathbb{R}[[\hbar]]} {^{\star}\tilde{\mathcal{T}}_U} \longrightarrow \mathcal{A}_U
\end{align*}
such that the matrix
\begin{align*}
({^{\star}g_{ij}}) \in \mathcal{A}_U^{n \times n}, \quad {^{\star}g_{ij}} = {^{\star}g}(E_i, \tilde{E}_j)
\end{align*}
is $\star$-invertible, i.e., there exists a unique matrix $({^{\star}g^{ij}}) \in \mathcal{A}_U ^{n \times n}$ such that
\begin{align*}
{^{\star}g_{ik}} \star {^{\star}g^{kj}} = {^{\star}g^{jk}} \star {^{\star}g_{ki}} = \delta_i^j.
\end{align*}
\end{defn}

By the associativity of star product $\star$, $({^{\star}g_{ij}})$ is invertible if and only if it has a left $\star$-inverse and a right $\star$-inverse. On the other hand, only associativity of the Moyal product is required in the proof of Proposition \ref{gij0}. Therefore, it holds true also replacing the Moyal product by the star product. Thus, on $U$,
\begin{align*}
({^{\star}g_{ij}})\,\, \mbox{is invertible} \Llrw ({^{\star}g_{ij}}[0])\,\, \mbox{is invertible}.
\end{align*}

In terms of noncommutative $\star$-metric $^{\star}g$, it can induce dual bases $E^i$, $\tilde{E}^j$ of $\tilde{E}_j$, $E_i$ respectively, which satisfy
\begin{align*}
^{\star}g(E^i, \tilde{E}_j) = {^{\star}g}(E_j, \tilde{E}^i)= \delta_j^i.
\end{align*}

\begin{defn}
The noncommutative left ({\rm{resp.}} right) $\star$-cotangent bundle $^{\star}\mathcal{T}_U^*$ ({\rm{resp.}} $^{\star}\tilde{\mathcal{T}}_U^*$) on $U$ with respect to the noncommutative metric $^{\star}g$ is the free left ({\rm{resp.}} right) $(\mathcal{A}_U, \star)$-module with basis $\{E^1, \cdots, E^n\}$ ({\rm{resp.}} $\{\tilde{E}^1, \cdots, \tilde{E}^n\}$)
\begin{align*}
^{\star}\mathcal{T}_U^* & = \Big\{a_i \star E^i \,\Big | \,a_i \in \mathcal{A}_U, \, a_i \star E^i = 0 \Longleftrightarrow a_i = 0 \Big\},\\
^{\star} \tilde{\mathcal{T}}_U^* & = \Big\{\tilde{E}^i \star \tilde{a}_i \,\Big | \, \tilde{a} _i \in \mathcal{A}_U, \, \tilde{E}^i \star \tilde{a}_i = 0 \Longleftrightarrow \tilde{a}_i = 0 \Big\}.
\end{align*}
\end{defn}

\begin{defn}
A noncommutative left ({\rm{resp.}} right) quasi-connection $^{\star} \nabla$ is a map
\begin{align*}
^{\star}\nabla: {^{\star}\mathcal{T}}_U \longrightarrow {^{\star}\tilde{\mathcal{T}}} ^* _U \otimes _{\mathcal{A}_U} {^{\star}\mathcal{T}}_U \quad
({\rm{resp.}}\,\,{^{\star}\tilde{\nabla}}: {^{\star}\tilde{\mathcal{T}}} _U \longrightarrow {^{\star}\tilde{\mathcal{T}}} _U \otimes _{\mathcal{A}_U} {^{\star}\mathcal{T}} ^* _U)
\end{align*}
such that noncommutative left ({\rm resp.} right) covariant derivatives
\begin{align*}
^{\star}\nabla_i : {^{\star}\mathcal{T}}_U \longrightarrow {^{\star}\mathcal{T}}_U \quad
({\rm{resp.}}\,\,{^{\star}\tilde{\nabla}}_i: {^{\star}\tilde{\mathcal{T}}} _U \longrightarrow {^{\star}\tilde{\mathcal{T}}} _U)
\end{align*}
defined by
\begin{align*}
^{\star}\nabla_i V = {^{\star}g}(E_i, \tilde{E} ^k) \star W_k \quad ({\rm{resp.}} \,\,{^{\star}\tilde{\nabla}}_i \tilde{V} =\tilde{W} _k \star {^{\star}g}(E^k, \tilde{E}_i))
\end{align*}
for any
\begin{align*}
{^{\star} \nabla} V =\tilde{E} ^k \otimes W_k,\,\, W_k \in {^{\star}\mathcal{T}}_U \quad  ({\rm{resp.}} \,\,{^{\star}\tilde{\nabla}} \tilde{V} =\tilde{W} _k \otimes E^k,\,\, \tilde{W}_k \in {^{\star}\tilde{\mathcal{T}}}^*_U)
\end{align*}
are $\mathbb{R}[[\hbar]]$-linear.
\end{defn}

\begin{defn}
The left ({\rm{resp.}} right) curvature operators for left ({\rm{resp.}} right) quasi-connections are defined as follows.
\begin{align*}
    ^{\star}\mathcal{R}_{E_iE_j} &= [^{\star}\nabla_i, {^{\star}\nabla}_j]:\,\, ^{\star}\mathcal{T}_U \longrightarrow {^{\star}\mathcal{T}}_U,\\
    ^{\star}\tilde{\mathcal{R}}_{\tilde{E}_i\tilde{E}_j} &= [^{\star}\tilde{\nabla}_i, {^{\star}\tilde{\nabla}}_j] : \,\,^{\star}\tilde{\mathcal{T}}_U \longrightarrow {^{\star}\tilde{\mathcal{T}}}_U.
\end{align*}
\end{defn}

Note that the $\star$-curvature operators are not $(\mathcal{A}_U, \star)$-linear because the $\star$-covariant derivatives do not satisfy the Leibniz rule. But the $\star$-Riemannian curvatures still can be defined formally.
\begin{defn}
The left ({\rm{resp.}} right) Riemannian curvatures for left ({\rm{resp.}} right) quasi-connections are defined as follows.
\begin{align*}
^{\star}R_{lkij} & = {^{\star}g}({^{\star}\mathcal{R}}_{E_iE_j}E_k, \tilde{E}_l), \\
^{\star}\tilde{R}_{lkij} & = -{^{\star}g}(E_k, {^{\star}\tilde{R}}_{\tilde{E}_i\tilde{E}_j}\tilde{E}_l).
\end{align*}
\end{defn}

Unlike the situation of the Moyal product, the quasi-connections are not compatible with $\star$-metrics, therefore two Riemannian curvatures with respect to star products are not equal in general. This yields four different noncommutative Ricci curvatures.
\begin{defn}
The left ({\rm{resp.}} right) Ricci curvatures for left ({\rm{resp.}} right) quasi-connections are defined as follows.
\begin{align*}
    ^{\star}R_{kj} = & {^{\star}g}({^{\star}\mathcal{R}}_{E_iE_j}E_k, \tilde{E}_l) \star {^{\star}g}^{li} = {^{\star} R}_{lkij} \star {^{\star}g}^{li},\\
    ^{\star}\Theta_{il} = &{^{\star}g}^{jk} \star {^{\star} g}({^{\star}\mathcal{R}}_{E_iE_j}E_k, \tilde{E}_l) = {^{\star}g}^{jk} \star {^{\star} R}_{lkij},\\
    ^{\star}\tilde{R}_{kj} = &- {^{\star}g}(E_k, {^{\star}\tilde{R}}_{\tilde{E}_i\tilde{E}_j}\tilde{E}_l) \star {^{\star}g}^{li} = {^{\star}\tilde{R}}_{lkij} \star {^{\star}g}^{li},\\
    ^{\star}\tilde{\Theta}_{il} = &- {^{\star}g}^{jk} \star {^{\star} g}(E_k, {^{\star}\tilde{R}}_{\tilde{E}_i\tilde{E}_j}\tilde{E}_l) = {^{\star}g}^{jk} \star {^{\star}\tilde{R}}_{lkij},\\
    ^{\star}R_j^p = &{^{\star}g}^{pk} \star {^{\star}g}({^{\star}\mathcal{R}}_{E_iE_j}E_k, \tilde{E}_l) \star {^{\star}g}^{li}= {^{\star}g}^{pk} \star {^{\star}R}_{lkij} \star {^{\star}g}^{li},\\
    ^{\star}\Theta_i^p = &{^{\star}g}^{jk} \star {^{\star}g}({^{\star}\mathcal{R}}_{E_iE_j}E_k, \tilde{E}_l) \star {^{\star}g}^{lp} = {^{\star}g}^{jk} \star {^{\star}R}_{lkij} \star {^{\star}g}^{lp},\\
    ^{\star}\tilde{R}^p_j = &- {^{\star}g}^{pk} \star {^{\star}g}(E_k, {^{\star}\tilde{R}}_{\tilde{E}_i\tilde{E}_j}\tilde{E}_l) \star {^{\star}g}^{li} = {^{\star}g}^{pk} \star {^{\star}\tilde{R}}_{lkij} \star {^{\star}g}^{li},\\
    ^{\star}\tilde{\Theta}^p_i = &- {^{\star}g}^{jk} \star {^{\star}g}(E_k, {^{\star}\tilde{R}}_{\tilde{E}_i\tilde{E}_j}\tilde{E}_l) \star {^{\star}g}^{lp} = {^{\star}g}^{jk} \star {^{\star}\tilde{R}}_{lkij} \star {^{\star}g}^{lp}.
\end{align*}
\end{defn}

The above definitions of Ricci curvatures have the same left and right traces. This yields the left and right scalar curvatures.
\begin{defn}
The left ({\rm{resp.}} right) scalar curvatures for left ({\rm{resp.}} right) quasi-connections are defined as follows.
\begin{align*}
^{\star}R = {^{\star}R}^j_j = {^{\star}\Theta}^i_i,\quad {^{\star}\tilde{R}} = {^{\star}\tilde{R}}^j_j = {^{\star}\tilde{\Theta}}^i_i.
\end{align*}
\end{defn}

\begin{prop}
Suppose that the left ({\rm{resp.}} right) quasi-connection is torsion free, i.e., it satisfies that
\begin{align*}
^{\star}\nabla_i E_j = {^{\star}\nabla}_j E_i,\quad
    ({\rm{resp.}} \,\,{^{\star} \tilde{\nabla}}_i \tilde{E}_j = {^{\star}\tilde{\nabla}}_j \tilde{E}_i).
\end{align*}
Then the first (algebraic) Bianchi identity
\begin{align*}
^{\star}\mathcal{R}_{E_iE_j}E_k + {^{\star}\mathcal{R}}_{E_jE_k}E_i + {^{\star}\mathcal{R}}_{E_kE_i}E_j &= 0,\\
{^{\star}\mathcal{\tilde{R}}}_{E_iE_j}E_k + {^{\star}\mathcal{\tilde{R}}}_{E_jE_k}E_i + {^{\star}\mathcal{\tilde{R}}}_{E_kE_i}E_j &= 0.
\end{align*}
holds for $1 \leq i, j, k \leq n$.
\end{prop}
\pf The proof of the first Bianchi identity in Theorem \ref{bianchi} can be applied to torsion free quasi-connections directly. \qed

Now we study quasi-connections given by an isometric embedding.
\begin{align*}
X: (U, g[0]) \longrightarrow \mathbb{R}^{p, m - p}.
\end{align*}
For $Y = (Y^1, \cdots, Y^m)$, $Z = (Z^1, \cdots, Z^m) \in \mathcal{A}_U^m$, we denote
\begin{align*}
Y \star_{\eta} Z = \sum\limits_{\alpha = 1}^m\eta_{\alpha \alpha} Y^{\alpha} \star Z^{\alpha}.
\end{align*}
Then the $\star$-metric is
\begin{align*}
    {^{\star}g}\big(E_i, \tilde{E}_j \big)=\partial_i X \star_{\eta} \partial _j X.
\end{align*}

\begin{lem}\label{sigma}
For isometric embedding $X$, the left ({\rm{resp.}} right) $(\mathcal{A}_U, \star)$-module homomorphism
\begin{align*}
\sigma : {^{\star}\mathcal{T}}_U \longrightarrow \mathcal{A}_U^m \quad
({\rm{resp.}} \,\,\tilde{\sigma} : {^{\star} \tilde{\mathcal{T}}}_U \longrightarrow \mathcal{A}_U^m)
\end{align*}
given by
\begin{align*}
\sigma(E_i) = \partial_i X, \quad ({\rm{resp.}} \,\,\tilde{\sigma}(\tilde{E}_i) =\partial_i X)
\end{align*}
is injective.
\end{lem}
\pf Assume
\begin{align*}
a^i \star \partial_i X = \sigma (a^i \star E_i) = 0.
\end{align*}
We obtain
\begin{align*}
0 = \sum\limits_{\alpha = 1}^m \eta _{\alpha\alpha} a^i \star \partial_i X^\alpha \star \partial _j X^\alpha = a^i \star {^{\star}g}_{ij}
\end{align*}
for each $j$. Hence, for each $k$,
\begin{align*}
0 = a^i \star {^{\star}g}_{ij} \star {^{\star}g}^{jk} = a^k
\end{align*}
Therefore $\sigma$ is injective. Same argument gives that $\tilde \sigma$ is injective. \qed

Denote
\begin{align*}
{^{\star}\mathcal{N}}_U & = \Big\{Y \in \mathcal{A}_U^m \Big| Y \star_{\eta} \tilde{\sigma}(\tilde{E}_i) = 0, \forall \,i\Big\},\\
{^{\star}\tilde{\mathcal{N}}}_U &= \Big\{Y \in \mathcal{A}_U^m \Big| \sigma(E_i)\star_{\eta} Y = 0, \forall \,i\Big\}.
\end{align*}
It is clear that ${^{\star}\mathcal{N}}_U$ (resp. ${^{\star}\tilde{\mathcal{N}}}_U$) is a left (resp. right) $(\mathcal{A}_U, \star)$-module.
\begin{lem}
For isometric embedding $X$, there is a direct sum decompositions of left ({\rm{resp.}} right) $(\mathcal{A}_U, \star)$-modules
    \begin{align}
        \mathcal{A}_U^m = \sigma({^{\star}\mathcal{T}}_U) \oplus {^{\star}\mathcal{N}}_U \quad
        ({\rm{resp.}}\,\,\mathcal{A}_U^m = \tilde{\sigma}({^{\star}\tilde{\mathcal{T}}}_U) \oplus {^{\star}\tilde{\mathcal{N}}}_U).
        \label{t-n-decom}
    \end{align}
\end{lem}
\pf Let $Y \in \mathcal{A}_U^m$. Denote
\begin{align*}
y^i = Y \star_{\eta} \tilde{\sigma}(\tilde{E}_j) \star {^{\star}g}^{ji}, \quad Y^{\top} = y^i \star \sigma(E_i), \quad Y^{\bot} = Y - Y^{\top}.
\end{align*}
Since, for each $j$,
\begin{align*}
    Y^{\bot} \star_{\eta} \tilde{\sigma}(\tilde{E}_j) = &Y \star_{\eta} \tilde{\sigma}(\tilde{E}_j) - Y^{\top} \star_{\eta} \tilde{\sigma}(\tilde{E}_j)\\
    = &Y \star_{\eta} \tilde{\sigma}(\tilde{E}_j) - y^i \star \sigma(E_i) \star_{\eta} \tilde{\sigma}(\tilde{E}_j)\\
    = &Y \star_{\eta} \tilde{\sigma}(\tilde{E}_j) - y^i \star {^{\star}g}_{ij}\\
    = &Y \star_{\eta} \tilde{\sigma}(\tilde{E}_j) - Y \star_{\eta} \tilde{\sigma}(\tilde{E}_k) \star {^{\star}g}^{ki} \star {^{\star}g}_{ij}\\
    = &Y \star_{\eta} \tilde{\sigma}(\tilde{E}_j) - Y \star_{\eta} \tilde{\sigma}(\tilde{E}_j)\\
    = &0,
\end{align*}
we find that
\begin{align*}
Y^{\bot} \in {^{\star}\mathcal{N}}_U.
\end{align*}
Since $({^{\star}g}_{ij})$ is invertible, we have
\begin{align*}
\sigma({^{\star}\mathcal{T}}_U) \cap {^{\star}\mathcal{N}}_U = \{0\},
\end{align*}
thus
\begin{align*}
\mathcal{A}_U^m = \sigma({^{\star}\mathcal{T}}_U) \oplus {^{\star}\mathcal{N}}_U.
\end{align*}
The case for the right modules can be proved by the same argument. \qed

Denote
\begin{align*}
{\rm pr}_1 : \mathcal{A}_U^m \longrightarrow \sigma({^{\star}\mathcal{T}}_U) \quad
({\rm{resp.}}\,\,{\rm \tilde{pr}}_1 : \mathcal{A}_U^m \longrightarrow \tilde{\sigma}({^{\star}\tilde{\mathcal{T}}}_U))
\end{align*}
the projection onto the first factor with respect to the decomposition (\ref{t-n-decom}). In terms of the projection, the left ({\rm{resp.}} right) quasi-connection
\begin{align*}
^{\star}\nabla: {^{\star}\mathcal{T}}_U \longrightarrow {^{\star}\tilde{\mathcal{T}}} ^* _U \otimes _{\mathcal{A}_U} {^{\star}\mathcal{T}}_U \quad
({\rm{resp.}}\,\,{^{\star}\tilde{\nabla}}: {^{\star}\tilde{\mathcal{T}}} _U \longrightarrow {^{\star}\tilde{\mathcal{T}}} _U \otimes _{\mathcal{A}_U} {^{\star}\mathcal{T}} ^* _U)
\end{align*}
is given by
\begin{align}
{^{\star}\nabla}V = \tilde{E}^k \otimes {^{\star}\nabla}_kV \quad ({\rm resp.}\,\, {^{\star}\tilde{\nabla}}\tilde{V} = {^{\star}\tilde{\nabla}}_k\tilde{V} \otimes E^k) \label{X-quasi-conn}
\end{align}
for any $V \in {^{\star}\mathcal{T}}_U \quad ({\rm resp.} \,\, \tilde{V} \in {^{\star}\tilde{\mathcal{T}}}_U)$, where the left ({\rm{resp.}} right) quasi-covariant derivative
\begin{align*}
^{\star}\nabla_i : {^{\star}\mathcal{T}}_U \longrightarrow {^{\star}\mathcal{T}}_U \quad
({\rm{resp.}}\,\,^{\star}\tilde{\nabla}_i : {^{\star}\tilde{\mathcal{T}}}_U \longrightarrow {^{\star}\tilde{\mathcal{T}}}_U),
\end{align*}
for each $i$, is given by
\begin{align}
^{\star}\nabla_i(V) = \sigma^{-1}\big({\rm pr}_1 (\partial_i\sigma(V))\big)\quad
\big({\rm{resp.}}\,\,^{\star}\tilde{\nabla}_i(\tilde{V}) = \tilde{\sigma}^{-1}\big({\rm \tilde{pr}}_1(\partial_i \tilde{\sigma}( \tilde{V}))\big)\big).\label{X-quasi-covar}
\end{align}

\begin{prop}
For isometric embedding $X$, the left ({\rm{resp.}} right) quasi-connection $^{\star}\nabla$ ({\rm{resp.}} $^{\star}\tilde{\nabla}$)
given by (\ref{X-quasi-covar}) and (\ref{X-quasi-conn}) is torsion free, i.e.,
    \begin{align*}
        ^{\star} \nabla_i E_j =  {^{\star} \nabla}_j E_i \quad
       ({\rm{resp.}}\,\, {^{\star}\tilde{\nabla}}_i \tilde{E}_j = {^{\star}\tilde{\nabla}}_j \tilde{E}_i).
    \end{align*}
\end{prop}
\pf The proposition follows from (\ref{X-quasi-covar}) and
\begin{align*}
{\rm pr}_1 \big(\partial_i \sigma(E_j)\big) & = {\rm pr}_1 (\partial_i \partial_j X),\\
{\rm \tilde{pr}}_1 \big(\partial_i \tilde{\sigma}(\tilde{E}_j)\big) & = {\rm \tilde{pr}}_1 (\partial_i \partial_j X).
\end{align*}
\qed

\appendix
\section{Moyal products of trigonometric functions}
\ls
Let $U \subset \mathbb{R}^2$ be an open subset with coordinates $(x^1, x ^2)=(\theta_1, \theta_2)$. Define the Moyal product $*$ on $\mathcal{A}_U = C^{\infty}(U)[[\hbar]]$ by the matrix
$$\begin{pmatrix}
0 & \lambda\\
- \lambda & 0
\end{pmatrix},$$
for some constant $\lambda \neq 0$. Moyal products of trigonometric functions are given as follows.
\subsection{}
\begin{align*}
    &(\sin\theta_1\sin\theta_2) * (\sin\theta_1\sin\theta_2)\\
    &\quad = \sin^2\theta_1\sin^2\theta_2\cosh^2(\lambda \hbar) - \cos^2\theta_1\cos^2\theta_2\sinh^2(\lambda \hbar),\\
    &(\sin\theta_1\sin\theta_2) * (\sin\theta_1\cos\theta_2)\\
    &\quad = \sin\theta_2\cos\theta_2\big(\sin^2\theta_1 + \sinh^2(\lambda \hbar)\big) - \sin\theta_1\cos\theta_1\cosh(\lambda \hbar)\sinh(\lambda \hbar),\\
    &(\sin\theta_1\sin\theta_2) * (\cos\theta_1\sin\theta_2)\\
    &\quad = \sin\theta_1\cos\theta_1\big(\sin^2\theta_2 + \sinh^2(\lambda \hbar)\big) + \sin\theta_2\cos\theta_2\cosh(\lambda \hbar)\sinh(\lambda \hbar),\\
    &(\sin\theta_1\sin\theta_2) * (\cos\theta_1\cos\theta_2)\\
    &\quad = \sin\theta_1\cos\theta_1\sin\theta_2\cos\theta_2 + (\sin^2\theta_1 - \sin^2\theta_2)\cosh(\lambda \hbar)\sinh(\lambda \hbar),
\end{align*}
\subsection{}
\begin{align*}
    &(\sin\theta_1\cos\theta_2) * (\sin\theta_1\sin\theta_2)\\
    &\quad = \sin\theta_2\cos\theta_2\big(\sin^2\theta_1 + \sinh^2(\lambda \hbar)\big) + \sin\theta_1\cos\theta_1\cosh(\lambda \hbar)\sinh(\lambda \hbar),\\
    &(\sin\theta_1\cos\theta_2) * (\sin\theta_1\cos\theta_2)\\
    &\quad = \sin^2\theta_1\cos^2\theta_2\cosh^2(\lambda \hbar) - \cos^2\theta_1\sin^2\theta_2\sinh^2(\lambda \hbar),\\
    &(\sin\theta_1\cos\theta_2) * (\cos\theta_1\sin\theta_2)\\
    &\quad = \sin\theta_1\cos\theta_1\sin\theta_2\cos\theta_2 + (\cos^2\theta_1 - \sin^2\theta_2)\cosh(\lambda \hbar)\sinh(\lambda \hbar),\\
    &(\sin\theta_1\cos\theta_2) * (\cos\theta_1\cos\theta_2)\\
    &\quad = \sin\theta_1\cos\theta_1\big(\cos^2\theta_2 + \sinh^2(\lambda \hbar)\big) - \sin\theta_2\cos\theta_2\cosh(\lambda \hbar)\sinh(\lambda \hbar),
\end{align*}
\subsection{}
\begin{align*}
    &(\cos\theta_1\sin\theta_2) * (\sin\theta_1\sin\theta_2)\\
    &\quad = \sin\theta_1\cos\theta_1\big(\sin^2\theta_2 + \sinh^2(\lambda \hbar)\big) - \sin\theta_2\cos\theta_2\cosh(\lambda \hbar)\sinh(\lambda \hbar),\\
    &(\cos\theta_1\sin\theta_2) * (\sin\theta_1\cos\theta_2)\\
    &\quad = \sin\theta_1\cos\theta_1\sin\theta_2\cos\theta_2 + (\sin^2\theta_1 - \cos^2\theta_2)\cosh(\lambda \hbar)\sinh(\lambda \hbar),\\
    &(\cos\theta_1\sin\theta_2) * (\cos\theta_1\sin\theta_2)\\
    &\quad = \cos^2\theta_1\sin^2\theta_2\cosh^2(\lambda \hbar) - \sin^2\theta_1\cos^2\theta_2\sinh^2(\lambda \hbar),\\
    &(\cos\theta_1\sin\theta_2) * (\cos\theta_1\cos\theta_2)\\
    &\quad = \sin\theta_2\cos\theta_2\big(\cos^2\theta_1 + \sinh^2(\lambda \hbar)\big) + \sin\theta_1\cos\theta_1\cosh(\lambda \hbar)\sinh(\lambda \hbar),
\end{align*}
\subsection{}
\begin{align*}
    &(\cos\theta_1\cos\theta_2) * (\sin\theta_1\sin\theta_2)\\
    &\quad = \sin\theta_1\cos\theta_1\sin\theta_2\cos\theta_2 + (\cos^2\theta_1 - \cos^2\theta_2)\cosh(\lambda \hbar)\sinh(\lambda \hbar),\\
    &(\cos\theta_1\cos\theta_2) * (\sin\theta_1\cos\theta_2)\\
    &\quad = \sin\theta_1\cos\theta_1\big(\cos^2\theta_2 + \sinh^2(\lambda \hbar)\big) + \sin\theta_2\cos\theta_2\cosh(\lambda \hbar)\sinh(\lambda \hbar),\\
    &(\cos\theta_1\cos\theta_2) * (\cos\theta_1\sin\theta_2)\\
    &\quad = \sin\theta_2\cos\theta_2\big(\cos^2\theta_1 + \sinh^2(\lambda \hbar)\big) - \sin\theta_1\cos\theta_1\cosh(\lambda \hbar)\sinh(\lambda \hbar),\\
    &(\cos\theta_1\cos\theta_2) * (\cos\theta_1\cos\theta_2)\\
    &\quad = \cos^2\theta_1\cos^2\theta_2\cosh^2(\lambda \hbar) - \sin^2\theta_1\sin^2\theta_2\sinh^2(\lambda \hbar).
\end{align*}

\bigskip

\footnotesize {\it Acknowledgement. The work is partially supported by the special foundation for Junwu and Guangxi Ba Gui Scholars.}

\end{document}